\documentclass{article}

\usepackage{microtype}
\usepackage{graphicx}
\usepackage{booktabs} 

\usepackage{hyperref}

\usepackage[accepted]{icml2024}

\usepackage{amsmath}
\usepackage{amssymb}
\usepackage{mathtools}
\usepackage{amsthm}

\usepackage[capitalize,noabbrev]{cleveref}

\theoremstyle{plain}
\newtheorem{theorem}{Theorem}[section]

\theoremstyle{definition}

\theoremstyle{remark}
\newtheorem{remark}[theorem]{Remark}
\newtheorem{property}[theorem]{Property}

\usepackage{bm}
\usepackage{multicol,multirow}
\usepackage{caption}
\usepackage{subcaption}
\newcommand\myrowlabel[1]{\rotatebox[origin=c]{90}{#1}}

\usepackage[textsize=tiny]{todonotes}

\icmltitlerunning{Verifying message-passing neural networks via topology-based bounds tightening}

\begin{document}

\twocolumn[
\icmltitle{Verifying message-passing neural networks via topology-based bounds tightening}



\icmlsetsymbol{equal}{*}

\begin{icmlauthorlist}
\icmlauthor{Christopher Hojny}{equal,tue}
\icmlauthor{Shiqiang Zhang}{equal,imp}
\icmlauthor{Juan S. Campos}{imp}
\icmlauthor{Ruth Misener}{imp}
\end{icmlauthorlist}

\icmlaffiliation{tue}{Eindhoven University of Technology, Eindhoven, The Netherlands}
\icmlaffiliation{imp}{Department of Computing, Imperial College London, UK}

\icmlcorrespondingauthor{Christopher Hojny}{c.hojny@tue.nl}

\icmlkeywords{Machine Learning, ICML}

\vskip 0.3in
]



\printAffiliationsAndNotice{\icmlEqualContribution} 

\begin{abstract}
    Since graph neural networks (GNNs) are often vulnerable to attack, we need to know when we can trust them. We develop a computationally effective approach towards providing robust certificates for message-passing neural networks (MPNNs) using a Rectified Linear Unit (ReLU) activation function. Because our work builds on mixed-integer optimization, it encodes a wide variety of subproblems, for example it admits (i) both adding and removing edges, (ii) both global and local budgets, and (iii) both topological perturbations and feature modifications. Our key technology, topology-based bounds tightening, uses graph structure to tighten bounds. We also experiment with aggressive bounds tightening to dynamically change the optimization constraints by tightening variable bounds. To demonstrate the effectiveness of these strategies, we implement an extension to the open-source branch-and-cut solver SCIP. We test on both node and graph classification problems and consider topological attacks that both add and remove edges.
\end{abstract}

\section{Introduction}\label{sec:introduction}
    Graph neural networks (GNNs) may have incredible performance in graph-based tasks, but researchers also raise concerns about their vulnerability: small input changes sometimes lead to wrong GNN predictions \cite{Gunnemann2022}. To study these GNN vulnerabilities, prior works roughly divide into two classes: adversarial attacks \cite{Dai2018,Zugner2018,Takahashi2019,Xu2019,Zugner2019b,Chen2020,Ma2020,Sun2020,Wang2020,Geisler2021} and certifiable robustness \cite{Bojchevski2019,Zugner2019a,Zugner2020,Bojchevski2020,Jin2020,Salzer2023}. Beyond the difficulties of developing adversarial robustness for dense neural networks \cite{Lomuscio2017,Fischetti2018}, incorporating graphs brings new challenges to both adversarial attacks and certifiable robustness. The first difficulty is defining graph perturbations because, beyond tuning the features \cite{Takahashi2019,Zugner2018,Zugner2019a,Bojchevski2020,Ma2020}, an attacker may inject nodes \cite{Sun2020,Wang2020} or add/delete edges \cite{Dai2018,Zugner2018,Bojchevski2019,Zugner2019b,Zugner2020,Xu2019,Chen2020,Jin2020,Geisler2021}. Second, binary elements in the adjacency matrix create discrete optimization problems. Finally, perturbations to a node may indirectly attack other nodes via message passing or graph convolution.

    \begin{figure*}[tp]
        \centering
        \includegraphics[width=.55\textwidth]{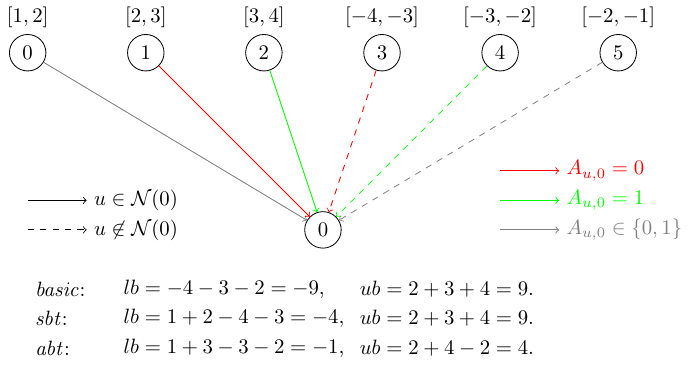}
        \includegraphics[width=.38\textwidth]{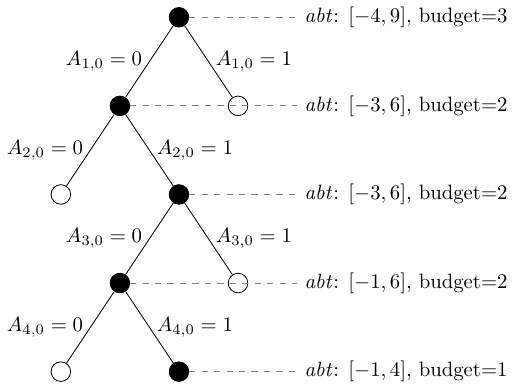}
        \vspace{-3mm}
        \caption{(\textbf{left}) Consider a graph with $6$ nodes $u = 0, \ldots, 5$ and one feature. The neighbor set of node $0$ is $\mathcal N(0)=\{0,1,2\}$. The input bounds are given above each node. Assume the budget, i.e., maximal number of modifications, for node $0$ is $3$. The modifications could be removing neighbors from $\mathcal N(0)$ or adding new neighbors from $\{3,4,5\}$. Four decisions have been made in the branch-and-bound tree, i.e., binary variables representing edges are set as $A_{1,0}=0,A_{2,0}=1,A_{3,0}=0,A_{4,0}=1$. Since node $2$ is a neighbor of node $0$ while node $3$ is not, fixing $A_{2,0}=1$ and $A_{3,0}=0$ spends no budget. For each method, we compute the bounds for node $0$ in the next layer. To compute a lower bound, the plain strategy (\textit{basic}) chooses all negative lower bounds without considering either budgets or previous decisions in the branch-and-bound tree. Static bounds tightening (\textit{sbt}), the first topology-based bounds tightening routine, removes node $2$ and adds node $3, 4$ as neighbors within $3$ budgets, but ignores decisions in the branch-and-bound tree. Aggressive bounds tightening (\textit{abt}) yields tighter bounds by saving node $0$ and adding node $5$ as neighbors. (\textbf{right}) The branch-and-bound tree corresponding to the left. We provide the bounds yielded from \textit{abt} and budget left after each decision.}
        \label{fig:overall}
    \end{figure*}

    Adversarial attacks aim to change the predictions of a GNN with admissible perturbations. In graph classification, the attacking goal is the prediction of a target graph \cite{Dai2018,Chen2020}. In node classification, attacks may be \emph{local} (or targeted) and \emph{global} (or untargeted). Local attacks \cite{Dai2018,Zugner2018,Takahashi2019,Wang2020} try to change the prediction of a single node under perturbations, and global attacks \cite{Xu2019,Zugner2019b,Ma2020,Sun2020,Geisler2021} allow perturbations to a group of nodes. Except for the Q-learning approach and genetic algorithm of \citet{Dai2018}, most aforementioned works are first-order methods which derive or approximate gradients w.r.t.\ features and edges. Binary variables are flipped when they are chosen to be updated.

    Certifiable robustness tries to guarantee that the prediction will not change under any admissible GNN perturbation. The state-of-the-art \cite{Bojchevski2019,Zugner2019a,Zugner2020,Jin2020} typically formulates certifiable robustness as a constrained optimization problem, where the objective is the worst-case margin between the correct class and other class(es), and the constraints represent admissible perturbations \cite{Gunnemann2022}. Given a GNN and a target node/graph, a certificate requires proving that the objective is always positive. Any feasible solution with a negative objective is an adversarial attack. Most existing certificates \cite{Zugner2019a,Jin2020} focus on graph convolutional networks (GCNs) \cite{Kipf2017}. The certificate on personalized propagation of neural predictions \cite{Gasteiger2019} relies on local budget and yields looser guarantees in the presence of a global budget \cite{Bojchevski2019}. Also, each certificate has specific requirements on the perturbation type, e.g., changing node features only \cite{Zugner2019a}, modifying graph structure only \cite{Bojchevski2019,Jin2020}, removing edges only \cite{Zugner2020}, and allowing only orthogonal Gromov-Wasserstein threats \cite{Jin2022}. Instead of verifying GNNs directly, several works \cite{Bojchevski2020,Wang2021,Xia2024} provide certified defenses based on the randomized smoothing framework \cite{Cohen2019}.

    This work develops certificates on the classic message passing framework, especially GraphSAGE \cite{Hamilton2017}. Using a recently-proposed mixed-integer formulation for GNNs \cite{McDonald2023,Zhang2023b,Zhang2023a}, we directly encode a GNN into an optimization problem using linear constraints. Many perturbations are compatible with our formulation:  (i) both adding and removing edges, (ii) both global and local budgets, and (iii) both topological perturbations and feature modifications. 

    When verifying fully-dense, feed-forward neural networks with ReLU activation, prior work shows that tightening variable bounds in a big-M formulation \cite{Anderson2020} may lead to better computational performance \cite{Tjeng2019,Botoeva2020,Tsay2021,Badilla2023,Zhao2024}. Since tighter variable bounds may improve the objective value of relaxations of the big-M formulation, they may be useful when providing a certificate of robustness. Because the optimization problems associated with verifying MPNNs are so large, this work cannot use the tighter, convex-hull based optimization formulations \cite{Singh2019a,Tjandraatmadja2020,Muller2022}. Instead, we use what we call \emph{topology-based bounds tightening} to enable a much stronger version of feasibility-based bounds tightening: we extend SCIP \cite{Bestuzheva2023} to explicitly use the graph structures. We also develop an \emph{aggressive bounds tightening} \cite{Belotti2016} routine to dynamically change the optimization constraints by tightening variable bounds within SCIP. Key outcomes include: (i) solving literature node classification instances in a fraction of a second, (ii) solving an extra 266 graph classification instances after implementing topology-based bounds tightening in SCIP, and (iii) making the open-source solver SCIP nearly as performant as the commercial solver Gurobi, e.g., improving the time penalty of the open-source solver from a factor of 10 to a factor of 3 for robust instances. 
    
    The paper begins in Section \ref{sec:GNN_encoding} by defining a mixed-integer encoding for MPNNs. Section \ref{sec:verification} presents the verification problem and develops our two topology-based bounds tightening routines, static bounds tightening \emph{sbt} and aggressive bounds tightening \emph{abt}. Section \ref{sec:experiments} presents the numerical experiments and Section \ref{sec:conclusion} concludes. Figure \ref{fig:overall} represents a toy example showing the basic approach in comparison to our two topology-based approaches \emph{sbt} and \emph{abt}.

\section{Definition \& Encoding of MPNNs}\label{sec:GNN_encoding}
    We inherit the mixed-integer formulation of MPNNs from \citet{Zhang2023b} and formulate ReLU activation using big-M \cite{Anderson2020}. Consider a trained GNN:
    \begin{equation}
        \begin{aligned}
            f:\mathbb R^{N\times d_0}\times\{0,1\}^{N\times N}&\to \mathbb R^{N\times d_L}\\
            (X,A)&\mapsto f(X,A)
        \end{aligned}
    \end{equation}
    whose $l$-th layer with weights ${\bm w}_{u\to v}^{(l)}$ and biases ${\bm b}_v^{(l)}$ is:
    \begin{align}\label{eq:fixed l-th GNN layer}
        {\bm x}_v^{(l)}=\sigma\left(\sum\limits_{u\in\mathcal N(v)\cup\{v\}}{\bm w}_{u\to v}^{(l)}{\bm x}_{u}^{(l-1)}+{\bm b}_v^{(l)}\right)
    \end{align}
    where $v\in V$, $V=\{0,1,\dots,N-1\}$ is the node set, $\mathcal N(v)$ is the neighbor set of node $v$, and $\sigma$ is activation. Given input features $\{{\bm x}_v^{(0)}\}_{v\in V}$ and the graph structure, we can derive the hidden features $\{{\bm x}_v^{(l)}\}_{v\in V},{\bm x}_v^{(l)} \in \mathbb R^{d_l}$. When $l=L$, we obtain the node representation ${\bm x}_v^{(L)}$ of each node. 

\subsection{Big-M formulation for MPNNs}\label{subsec:big-M of GNN}
    When the graph structure is not fixed, we need to include all possible contributions from all nodes, i.e., the $l$-th layer is:
    \begin{align}\label{eq:non-fixed l-th GNN layer}
        {\bm x}_v^{(l)}=\sigma\left(\sum\limits_{u\in V}A_{u,v}{\bm w}_{u\to v}^{(l)}{\bm x}_{u}^{(l-1)}+{\bm b}_v^{(l)}\right)
    \end{align}
    where $A_{u,v}\in\{0,1\}$ controls the existence of edge $u\to v$. 

    With fixed weights and biases, we still need to handle the nonlinearities caused by (i) bilinear terms $A_{u,v}{\bm x}_{u}^{(l-1)}$, and (ii) activation $\sigma$. Let $\bar {\bm x}_v^{(l)}=\sum\limits_{u\in V}A_{u,v}{\bm w}_{u\to v}^{(l)}{\bm x}_{u}^{(l-1)}+{\bm b}_v^{(l)}$, the \citet{Zhang2023a} big-M formulation introduces auxiliary variables ${\bm x}_{u\to v}^{(l-1)}$ to replace the bilinear terms $A_{u,v}{\bm x}_u^{(l-1)}$ and linearly encodes $\bar {\bm x}_v^{(l)}$:
    \begin{align}\label{eq:linear l-th GNN layer}
        \bar{\bm x}_v^{(l)}=\sum\limits_{u\in V}{\bm w}_{u\to v}^{(l)}{\bm x}_{u\to v}^{(l-1)}+{\bm b}_v^{(l)}.
    \end{align}
    Let $F_l:=\{0,1,\dots,d_l-1\}$ and denote the $f$-th element of ${\bm x}_{*}^{(l)}$ by $x_{*,f}^{(l)},f\in F_l$. Use $lb(\cdot)$ and $ub(\cdot)$ to represent the lower and upper bound of a variable, respectively. Then $x_{u\to v,f}^{(l-1)}=A_{u,v}x_{u,f}^{(l-1)}$ is equivalently formulated in the following big-M constraints:
    \begin{subequations}\label{eq:big-M of GNN}
        \begin{align}
            x_{u\to v,f}^{(l-1)}&\ge lb(x_{u,f}^{(l-1)})\cdot A_{u,v}\\
            x_{u\to v,f}^{(l-1)}&\le ub(x_{u,f}^{(l-1)})\cdot A_{u,v}\\
            x_{u\to v,f}^{(l-1)}&\le x_{u,f}^{(l-1)}-lb(x_{u,f}^{(l-1)})\cdot(1-A_{u,v})\\
            x_{u\to v,f}^{(l-1)}&\ge x_{u,f}^{(l-1)}-ub(x_{u,f}^{(l-1)})\cdot(1-A_{u,v}).
        \end{align}
    \end{subequations}

\subsection{Big-M formulation for ReLU}\label{subsec:big-M of ReLU}
    When using ReLU as the activation, i.e., 
    \begin{align}\label{eq:ReLU}
        x_{v,f}^{(l)}=\max\{0,\bar x_{v,f}^{(l)}\},
    \end{align}
    \citet{Anderson2020} proposed a big-M formulation:
    \begin{subequations}\label{eq:big-M of ReLU}
        \begin{align}
            x_{v,f}^{(l)}&\ge 0\\
            x_{v,f}^{(l)}&\ge \bar x_{v,f}^{(l)}\\
            x_{v,f}^{(l)}&\le \bar x_{v,f}^{(l)}-lb(\bar x_{v,f}^{(l)})\cdot(1-\sigma_{v,f}^{(l)})\label{eq:big-M of ReLU3}\\
            x_{v,f}^{(l)}&\le ub(\bar x_{v,f}^{(l)})\cdot\sigma_{v,f}^{(l)}
        \end{align}
    \end{subequations}
    where $\sigma_{v,f}^{(l)}\in\{0,1\}$ controls the on/off of the activation:
    \begin{align}
        x_{v,f}^{(l)}=
        \begin{cases}
            0,&\sigma_{v,f}^{(l)}=0\\
            \bar x_{v,f}^{(l)},&\sigma_{v,f}^{(l)}=1.
        \end{cases}
    \end{align}

\subsection{Bounds propagation}\label{subsec:bounds propagation}
    Eqs.\ \eqref{eq:big-M of GNN} and \eqref{eq:big-M of ReLU} show the importance of variable bounds $lb(\cdot)$ and $ub(\cdot)$ or \emph{big-M parameters}. Given the input feature bounds, we define bounds for auxiliary variables $x_{u\to v,f}^{(l-1)}$ and post-activation variables $x_{v,f}^{(l)}$. Using $x_{u\to v,f}^{(l-1)}=A_{u,v}x_{u,f}^{(l-1)}$, the bounds of $x_{u\to v,f}^{(l-1)}$ are:
    \begin{subequations}\label{eq:bounds of GNN auxiliary}
        \begin{align}
            lb(x_{u\to v,f}^{(l-1)})&=\min\{0,lb(x_{u,f}^{(l-1)})\}\\
            ub(x_{u\to v,f}^{(l-1)})&=\max\{0,ub(x_{u,f}^{(l-1)})\}
        \end{align}
    \end{subequations}
    with which we can use arithmetic propagation or feasibility-based bounds tightening to obtain bounds of $\bar x_{v,f}^{(l)}$ based on Eq.\ \eqref{eq:linear l-th GNN layer}. Then, we use Eq.\ \eqref{eq:ReLU} to bound $x_{v,f}^{(l)}$:
    \begin{subequations}\label{eq:bounds of ReLU auxiliary}
        \begin{align}
            lb(x_{v,f}^{(l)})&=\max\{0,lb(\bar x_{v,f}^{(l)})\}\\
            ub(x_{v,f}^{(l)})&=\max\{0,ub(\bar x_{v,f}^{(l)})\}.
        \end{align}
    \end{subequations}
    Without extra information, bounds defined in Eqs.\ \eqref{eq:bounds of GNN auxiliary} and \eqref{eq:bounds of ReLU auxiliary} are the tightest possible which derive from interval arithmetic. However, in a branch-and-bound tree, more and more variables will be fixed, which provides the opportunity to tighten the bounds. Additionally, in specific applications such as verification, the graph domain is restricted, allowing us to derive tighter bounds. 

    \begin{remark}
        A MPNN with $L$ message passing steps is suitable for node-level tasks. For graph-level tasks, there is usually a pooling layer after message passing to obtain a global representation and several dense layers thereafter as a final regressor/classifier. We omit these formulations since (i) linear pooling, e.g., mean and sum, is easily incorporated into our formulation, and (ii) dense layers are a special case of Eq.\ \eqref{eq:fixed l-th GNN layer} with a single node.
    \end{remark}

\section{Verification of MPNNs}\label{sec:verification}
\subsection{Problem definition}\label{subsec:problem definition}
    First, consider node classification. Given a trained MPNN defined as Eq.\ \eqref{eq:fixed l-th GNN layer}, the number of classes is the number of output features, i.e., $\mathcal C=d_L$, and the predicted label of node $t$ corresponds to the maximal logit, i.e., $c^*=\max_{c\in \mathcal C}f_{t,c}(X,A)$. Given an input $(X^*,A^*)$ consisting of features $X^*$ and adjacency matrix $A^*$, denote its predictive label for a target node $t$ as $c^*$. The worst case margin between predictive label $c^*$ and attack label $c$ under perturbations $\mathcal P(\cdot)$ is:
    \begin{equation}\label{eq:node classification}
        \begin{aligned}
            m^t(c^*,c):=\min\limits_{(X,A)}&~f_{t,c^*}(X,A)-f_{t,c}(X,A)\\
            s.t.&~ X\in\mathcal P(X^*),~A\in\mathcal P(A^*).
        \end{aligned}
    \end{equation}
    A positive $m^t(c^*,c)$ means that the logit of class $c^*$ is always larger than class $c$. If $m^t(c^*,c)>0,\forall c\in\mathcal C\backslash\{c^*\}$, then any admissible perturbation can not change the label assigned to node $t$, that is, this MPNN is robust to node $t$.

    For graph classification, instead of considering a single node, we want to know the worst case margin between two classes for a target graph, i.e.,
    \begin{equation}\label{eq:graph classification}
        \begin{aligned}
            m(c^*,c):=\min\limits_{(X,A)}&~f_{c^*}(X,A)-f_{c}(X,A)\\
            s.t.&~ X\in\mathcal P(X^*),~A\in\mathcal P(A^*).
        \end{aligned}
    \end{equation}
    In both problems, the target graph $(X^*,A^*)$ and predictive label $c^*$ are fixed. For node classification, the target node $t$ is also given. Therefore, we omit $t$ in Eq.\ \eqref{eq:node classification} and reduce both problems to one as shown in Eq.\ \eqref{eq:graph classification}.

\subsection{Admissible perturbations}\label{subsec:perturbations}
    The perturbations on features and edges can be described similarly. Locally, we may only change features/edges for each node with a given local budget. Also, there is typically a global budget for the number of changes. The feature perturbations are typically easier to implement since they will not hurt the message passing scheme, i.e., the graph structure is fixed. In such settings, there is no need to use the mixed-integer formulations for MPNNs in Section \ref{subsec:big-M of GNN} since a message passing step is actually simplified as a dense layer. Since feature perturbations are well-studied \cite{Zugner2019a} and our proposed bounds tightening techniques mainly focus on changeable graph structures, our computational results only consider the perturbations on the adjacency matrix.

    We first define the admissible perturbations for undirected graphs, which admits both adding and removing edges. Denote the global budget by $Q$ and local budget to node $v$ by $q_v$, then the perturbations $\mathcal P_1(A^*)$ are defined as:
    \begin{equation}\label{eq:undirected perturbations}
        \begin{aligned}
            \mathcal P_1(A^*)=\{&A\in\{0,1\}^{N\times N}~|~A=A^T,\\
            &\|A-A^*\|_0\le 2Q,\\
            &\|A_v-A_v^*\|_0\le q_v,~\forall v\in V\}
        \end{aligned}
    \end{equation}
    where $A_v$ is the $v$-th column of $A$. $\mathcal P_1(\cdot)$ will be used in graph classification since the graphs in benchmarks are usually undirected and relatively small.

    For node classification, literature benchmarks usually (i) are large directed graphs, e.g., $3000$ nodes, (ii) have many node features, e.g., $3000$ features, (iii) have small average degree, e.g., $1\sim 3$. If admitting adding edges, then the graph domain is too large to optimize over. Therefore, the state-of-the-art \cite{Zugner2020} only considers removing edges, where a $L$-hop neighborhood around the target node $t$ is sufficient. For a MPNN with $L$ message passing steps without perturbations, nodes outside a $L$-hop neighborhood cannot affect the prediction of $t$. Since adding edges is not allowed, the $L$-hop neighborhood of $t$ after perturbations is always a subset of the unperturbed neighborhood. Similar to the literature, we define a more restrictive perturbation space for large graphs:
    \begin{equation}\label{eq:directed perturbations}
        \begin{aligned}
            \mathcal P_2(A^*)=\{&A\in\{0,1\}^{N\times N}~|~\\
            &A_{u,v}\le A_{u,v}^*,~\forall u,v\in V,\\
            &\|A-A^*\|_0\le Q,\\
            &\|A_v-A_v^*\|_0\le q_v,~\forall v\in V\}
        \end{aligned}
    \end{equation}
    where global budget is replaced by $Q$ since the graph is directed. $\mathcal P_2(\cdot)$ will be used in node classification. For our later analysis, however, we focus on $\mathcal P_1(\cdot)$ since $\mathcal P_2(\cdot)$ is more like a special case without adding edges.

\subsection{Static bounds tightening}\label{subsec:sbt}

    \begin{algorithm}[tb]
       \caption{Static bounds tightening (\emph{sbt})}
       \label{alg:sbt}
        \begin{algorithmic}
           \STATE {\bfseries Input:} Input features ${\bm x}_v^{(0)}$, weights ${\bm w}_{u\to v}^{(l)}$, biases ${\bm b}_v^{(l)}$.
           \STATE Initialize $lb({\bm x}^{(0)})=ub({\bm x}^{(0)})={\bm x}^{(0)}$.
           \FOR{$l=1$ {\bfseries to} $L$}
                \STATE Get $lb(\bar {\bm x}^{(l)})$ using Eq.~\eqref{eq:decomposed definition of lower bound}.
                \COMMENT{\emph{basic} uses Eq.~\eqref{eq:basic lower bound}.}
                \STATE Get $lb({\bm x}^{(l)})$ using Eq.~\eqref{eq:bounds of ReLU auxiliary}.
                \STATE Get $ub(\bar{\bm x}^{(l)})$ and $ub({\bm x}^{(l)})$ in a similar way.
           \ENDFOR
        \end{algorithmic}
    \end{algorithm}

    Note that large budgets in Eq.\ \eqref{eq:undirected perturbations} make the verification problems meaningless since the perturbed graph could be any graph. The very basic assumption is that the perturbed graph is similar to the original one, which brings us to propose the first bounds tightening approach. The rough idea is to consider budgets when computing bounds of $\bar {\bm x}_{v}^{(l)}$ based on bounds of ${\bm x}_{u\to v}^{(l-1)}$ in Eq.\ \eqref{eq:linear l-th GNN layer}. Instead of considering all contributions from all nodes, we first calculate the bounds based on original neighbors and then maximally perturb the bounds with given budgets. Mathematically, $lb(\bar x_{v,f'}^{(l)})$ is found by solving the following optimization problem:
    \begin{equation}\label{eq:full definition of lower bound}
        \begin{aligned}
            \min\limits_{A,{\bm x}^{(l-1)}}~&\sum\limits_{u\in V}A_{u,v}\sum\limits_{f\in F_{l-1}}w_{u\to v,f\to f'}^{(l)}x_{u,f}^{(l-1)}+b_{v,f'}^{(l)}\\
            s.t.~& A\in \mathcal P_1(A^*)\\
            ~& {\bm x}^{(l-1)}\in [lb({\bm x}^{(l-1)}), ub({\bm x}^{(l-1)})]
        \end{aligned}
    \end{equation}
    where ${\bm x^{(l-1)}}:=\{x_{u,f}^{(l-1)}\}_{u\in V, f\in F_{l-1}}$, $w_{u\to v,f\to f'}^{(l)}$ is the $(f,f')$-th element in ${\bm w}_{u\to v}^{(l)}$.

    \begin{remark}
        For brevity, we omit the superscripts of layers for all variables, and subscripts of edges in weights, i.e., rewriting Eq.\ \eqref{eq:full definition of lower bound} as:
        \begin{equation}\label{eq:simplified definition of lower bound}
            \begin{aligned}
                lb(\bar x_{v,f'})=\min\limits_{A,{\bm x}}~&\sum\limits_{u\in V}A_{u,v}\sum\limits_{f\in F_{l-1}}w_{f,f'}x_{u,f}+b_{v,f'}\\
                s.t.~& A\in \mathcal P_1(A^*)\\
                ~& {\bm x}\in [lb({\bm x}), ub({\bm x})].
            \end{aligned}
        \end{equation}
    \end{remark}

    \begin{property}\label{prop:sbt}
        Eq.\ \eqref{eq:simplified definition of lower bound} is equivalent to:
        \begin{equation}\label{eq:decomposed definition of lower bound}
            \begin{aligned}
                lb(\bar x_{v,f'})=\sum\limits_{u\in \mathcal N^*(v)}lb_{u\to v}+b_{v,f'}+\min\limits_{|V_{lb}|\le q_v}\sum\limits_{u\in V_{lb}}\Delta_{u\to v}
            \end{aligned}
        \end{equation}
        where $\mathcal N^*(v)$ denotes the original neighbor set of node $v$, $lb_{u\to v}$ represents the contribution of node $u$ to the lower bound of node $v$ when $u$ is a neighbor of v, $\Delta_{u\to v}$ denotes the change of lower bound of node $v$ caused by modifying edge $u\to v$, $V_{lb}$ is the set of nodes consisting of removed/added neighbors of node $v$.
    \end{property}
    
    Calculating $lb_{u\to v}$ is straightforward:
    \begin{equation}\label{eq:contribution of neighbors}
        \begin{aligned}
            lb_{u\to v}=&\sum\limits_{f\in F_{l-1}}w_{f,f'}\cdot \mathbb I_{w_{f,f'}\ge 0}\cdot lb(x_{u,f})\\
            &+\sum\limits_{f\in F_{l-1}}w_{f,f'}\cdot \mathbb I_{w_{f,f'}< 0}\cdot ub(x_{u,f})
        \end{aligned}
    \end{equation}
    which is used to derive $\Delta_{u\to v}$ as:
    \begin{align}\label{eq:changes of add/remove neighbors}
        \Delta_{u\to v}=
        \begin{cases}
            -lb_{u\to v},&u\in \mathcal N^*(v)\\
            lb_{u\to v},&u\notin\mathcal N^*(v)
        \end{cases}
    \end{align}
    where two cases correspond to removing neighbor $u$ and adding $u$ as a neighbor, respectively. Furthermore, the last minimal term in Eq.\ \eqref{eq:decomposed definition of lower bound} is equivalent to choosing at most $q_v$ smallest negative terms among $\{\Delta_{u\to v}\}_{u\in V}$. The time complexity to compute lower bounds following Eq.\ \eqref{eq:decomposed definition of lower bound} for each feature in $l$-th layer is $O(N^2d_{l-1}d_l+N\log N)$. Upper bounds could be defined similarly, which are not included here due to space limitation.
    
    \begin{remark}
        The plain strategy without considering graph structure and budgets \textit{basic} is:
        \begin{equation}\label{eq:basic lower bound}
            \begin{aligned}
                lb(\bar x_{v,f'})=\sum\limits_{u\in V}\min\{0,lb_{u\to v}\}
            \end{aligned}
        \end{equation}
        and the time complexity is $O(N^2d_{l-1}d_l)$.
    \end{remark}
    
    Algorithm \ref{alg:sbt} calculates \emph{sbt} bounds (and \emph{basic} bounds) in a single forward pass of the model. As shown in the example presented in the Figure \ref{fig:overall} example, the bounds derived from \textit{basic} is $[-9,9]$, which is improved to $[-4,9]$ after applying static bounds tightening \textit{sbt}. 

    \begin{algorithm}[tb]
        \caption{Aggressive bounds tightening (\emph{abt})}
        \label{alg:abt}
        \begin{algorithmic}
            \STATE {\bfseries Input:} Input features ${\bm x}_v^{(0)}$, weights ${\bm w}_{u\to v}^{(l)}$, biases ${\bm b}_v^{(l)}$.
            \STATE Initialize $lb({\bm x}^{(0)})=ub({\bm x}^{(0)})={\bm x}^{(0)}$.
            \FOR{each node in branch-and-bound tree}
                \FOR{$l=1$ {\bfseries to} $L$}
                    \STATE Update $q'_v$ using Eq.~\eqref{eq:budgets in B&B}.
                    \STATE Update $lb(\bar {\bm x}^{(l)})$ using Eq.~\eqref{eq:decomposed lower bound in B&B}.
                    \STATE Update $lb({\bm x}^{(l)})$ using Eq.~\eqref{eq:bounds of ReLU auxiliary}.
                    \STATE Update $ub(\bar{\bm x}^{(l)})$ and $ub({\bm x}^{(l)})$ in a similar way.
               \ENDFOR
            \ENDFOR
        \end{algorithmic}
    \end{algorithm}

\subsection{Aggressive bounds tightening}\label{subsec:abt}
    Consider any node in a branch-and-bound tree, values of several $A_{u,v}$ are already decided during the path from root to current node, with which we can further tighten bounds in the subtree rooted by this node. \citet{Belotti2016} refer to the idea of tightening bounds in the branch-and-bound tree as \emph{aggressive bounds tightening}. In MPNN verification, there are three types of $A_{u,v}$ in Eq.\ \eqref{eq:simplified definition of lower bound}: (i) $A_{u,v}$ is fixed to $0$, (ii) $A_{u,v}$ is fixed to $1$, and (iii) $A_{u,v}$ is not fixed yet. Denote $V_0=\{u\in V~|~A_{u,v}=0\}$ and $V_1=\{u\in V~|~A_{u,v}=1\}$, then Eq.\ \eqref{eq:simplified definition of lower bound} in the current node is restricted as:
    
    \begin{equation}\label{eq:lower bound in B&B}
        \begin{aligned}
            lb(\bar x_{v,f'})=\min\limits_{A}~&\sum\limits_{u\in V}A_{u,v}lb_{u\to v}+b_{v,f'}\\
            s.t.~& A\in \mathcal P_1(A^*)\\
            & A_{u,v}=0,~\forall u\in V_0\\
            & A_{u,v}=1,~\forall u\in V_1.
        \end{aligned}
    \end{equation}
    
    \begin{property}\label{prop:abt}
        Eq.\ \eqref{eq:lower bound in B&B} is equivalent to:
        \begin{equation}\label{eq:decomposed lower bound in B&B}
            \begin{aligned}
                lb(\bar x_{v,f'})=&\sum\limits_{u\in (\mathcal N^*(v)\backslash V_0)\cup V_1}lb_{u\to v}+b_{v,f'}\\
                +&\min\limits_{V_{lb}\subseteq V\backslash(V_0\cup V_1), |V_{lb}|\le q_v'}\sum\limits_{u\in V_{lb}}\Delta_{u\to v}
            \end{aligned}
        \end{equation}
        where $q_v'$ is the currently available budget for node $v$.
    \end{property}
    
    In Eq.\ \eqref{eq:decomposed lower bound in B&B}, the first term sums over $(\mathcal N^*(v)\backslash V_0)\cup V_1$ to represent the contributions from current neighbors. The last term excludes all fixed edges and only considers changeable neighbors. Similarly, this minimal term equals to choose at most $q_v'$ smallest negative terms among $\{\Delta_{u\to v}\}_{u\in V\backslash(V_0\cup V_1))}$. $q_v'$ is derived from the remaining local and global budgets:
    \begin{equation}\label{eq:budgets in B&B}
        \begin{aligned}
            q_v'=\min\{&q_v-e_r(v)-e_a(v),\\
            &Q-\frac{1}{2}\sum\limits_{u\in V}e_r(u)-\frac{1}{2}\sum\limits_{u\in V}e_a(u)\}
        \end{aligned}
    \end{equation}
    where $e_r(v):=|\mathcal N^*(v)\cap V_0|$ is the number of removed edges around node $v$, $e_a(v):=|V_1\backslash \mathcal N^*(v)|$ is the number of added edges around node $v$.

    Algorithm \ref{alg:abt} describes the \emph{abt} steps. Since we derive \emph{abt} bounds within the branch-and-bound tree (when solvers will not directly change variable bounds), we add local cutting planes to implement \emph{abt} bounds (see Appendix \ref{app:implementation_details} example).

    \begin{remark}\label{remark:abt is sbt on each node}
        The bounds tightening inside the branch-and-bound tree can be interpreted as applying the Section \ref{subsec:sbt} bounds tightening to a modified target graph with a reduced budget. The spent budget changes the neighbors of node $v$ from $\mathcal N^*(v)$ to $(\mathcal N^*(v)\backslash V_0)\cup V_1$.
    \end{remark}

    As shown in Figure \ref{fig:overall}, aggressive bounds tightening \textit{abt} gives tighter bounds $[-1,4]$ comparing to \textit{basic} and \textit{sbt}. The branch-and-bound tree in Figure \ref{fig:overall} shows \textit{abt} in action.

   \begin{table}[h]
        \caption{Information on benchmarks. For multiple graphs, we compute the average number of nodes and edges.}
        \label{tab:benchmarks}
        \vskip 0.15in
        \begin{center}
        \begin{small}
        \setlength{\tabcolsep}{.45em}
        \begin{tabular}{lrrrrr}
            \toprule
            benchmark & \#graphs & \#nodes & \#edges & \#features & \#classes\\
            \midrule
            MUTAG & 188 & $17.9$ & $39.6$ & $7$ & $2$ \\
            ENZYMES & $600$ & $32.6$ & $124.3$ & $3$ & $6$ \\
            \midrule
            Cora & $1$ & $2708$ & $5429$ & $1433$ & $7$ \\
            CiteSeer & $1$ & $3312$ & $4715$ & $3703$ & $6$\\
            \bottomrule
        \end{tabular}
        \end{small}
        \end{center}
        \vskip -0.1in
    \end{table}

\subsection{Bounds tightening for ReLU}
Aggressive bounds tightening could also tighten ReLU interval bounds in feed-forward neural networks (NNs). Although we did not directly tighten ReLU bounds in MPNNs, local cutting planes representing tighter big-M coefficients $lb(\bar x_{v,f}^{(l)}), ub(\bar x_{v,f}^{(l)})$ in Eq.~\eqref{eq:big-M of ReLU} are added after applying \emph{abt}. For each ReLU, if its pre-activation variable $x_{v,f}^{(l)}$ has a non-negative lower bound or a non-positive upper bound, the binary variable $\sigma_{v,f}^{(l)}$ controlling on/off of this ReLU will be fixed due to these local cutting planes. Therefore, \emph{abt} implies a dynamic tightening on ReLU interval bounds in MPNNs. This idea could be applied to NNs since a NN is an MPNN with a single node in the graph.

For ReLU NNs, bounds tightening techniques yielding tighter bounds than interval arithmetic include: FastLin \cite{Weng2018}, CROWN \cite{Zhang2018}, DeepPoly \cite{Singh2018, Singh2019b}, and optimization-based bound tightening (OBBT) \cite{Tjeng2019,Tsay2021}. But these techniques are rarely applied in GNN verification. Although a few works \cite{Zugner2019a, Jin2020} involve convex relaxations for ReLUs, they do not tighten bounds for ReLUs.

OBBT is the only one of the existing advanced approaches which could immediately apply to GNNs. OBBT yields tighter bounds with high computational cost of solving many linear programs (LPs) or mixed-integer programs (MIPs). Incorporating computationally-effective  methods like FastLin, CROWN or DeepPoly is difficult for GNNs because we lose the linearity between layers. Such linearity is crucial to these approaches, e.g., the linear lower/upper bounds in FastLin/CROWN, or the zonotopy abstraction of DeepPoly. When the input graph structure is fixed, bounds tightening on ReLUs in GNNs is equivalent to its counterpart in NNs \citep{Wu2022}. But for non-fixed graphs, the links between two adjacent layers in GNNs are controlled by an adjacency matrix, whose elements are binary variables. So it is considering topological perturbations that loses the linearity between layers and inspires \emph{sbt} and \emph{abt}.

\section{Experiments}\label{sec:experiments}

    \begin{table*}[t]
        \caption{Summary of results for graph classification with local attack strength $s=2$. The approaches tested are SCIP basic (SCIPbasic), SCIP static bounds tightening (SCIPsbt), SCIP aggressive bounds tightening (SCIPabt), Gurobi basic (GRBbasic), and Gurobi static bounds tightening (GRBsbt). For each global budget, the number of instances is 600 for ENZYMES and 188 for MUTAG, but we only present comparisons when all methods give consistent results except for time out, as shown in column ``\#". We compare times with respect to both average (``avg-time") and the shifted geometric mean (``sgm-time"), as well as the number of solved instances within time limits 2 hours (``\# solved"). Since robust instances are the ones where mixed-integer performance is most important (non-robust instances may frequently be found by more heuristic approaches), we also compare on the set of robust instances.}
        \label{tab:summary_graph_classification_s2}
        \begin{center}
        \begin{small}
        \begin{tabular*}{\textwidth}{@{\;\;\extracolsep{\fill}}llrrrrrrrr}
            \toprule
            \multirow{2}{*}{benchmark} & \multirow{2}{*}{method} & \multicolumn{4}{c}{all instances} & \multicolumn{4}{c}{robust instances}\\
            \cmidrule{3-6} \cmidrule{7-10}
             & & \# & avg-time(s) & sgm-time(s) & \# solved  & \# & avg-time(s) & sgm-time(s) & \# solved\\
             
            \midrule
            \multirow{5}{*}{ENZYMES}
            & SCIPbasic  & 5915 & 605.97 & 37.81 & 5579 & 3549 & 278.58 & 12.51 & 3444 \\
            & SCIPsbt    & 5915 & \textbf{230.59} & \textbf{21.26} & \textbf{5831} & 3549 &  \textbf{82.89} & 
 \textbf{6.65} & \textbf{3528} \\
            & SCIPabt    & 5915 & 246.02 & 21.77 & 5817 & 3549 &  88.95 &  6.71 & 3522 \\
            \cmidrule{2-10}
            & GRBbasic   & 5915 &  86.03 &  7.59 & 5892 & 3549 &  32.80 &  2.82 & 3542 \\
            & GRBsbt     & 5915 &  \textbf{74.87} &  \textbf{7.09} & \textbf{5895} & 3549 &  \textbf{22.90} &  \textbf{2.50} & \textbf{3548} \\
            \midrule
            \multirow{5}{*}{MUTAG}
            & SCIPbasic  & 1589 & 679.86 & 189.75 & 1575 & 44 & 798.47 & 202.93 & 40 \\
            & SCIPsbt    & 1589 & \textbf{196.07} &  \textbf{75.17} & \textbf{1589} & 44 & 336.41 & 100.86 & \textbf{44} \\
            & SCIPabt    & 1589 & 207.50 &  82.43 & \textbf{1589} & 44 & \textbf{238.10} &  \textbf{91.11} & \textbf{44} \\
            \cmidrule{2-10}
            & GRBbasic   & 1589 &  \textbf{34.58} &  \textbf{4.06} & \textbf{1589} & 44 & 162.25 &  \textbf{12.11} & \textbf{44} \\
            & GRBsbt     & 1589 &  59.93 & 22.40 & 1589 & 44 &  \textbf{73.78} &  15.00 &  \textbf{44} \\
            \bottomrule
        \end{tabular*}
        \end{small}
        \end{center}
        \vskip -0.1in
    \end{table*}

     \begin{figure*}[ht]
        \centering
        \begin{subfigure}[b]{0.32\textwidth}
            \centering
            \includegraphics[width=\textwidth]{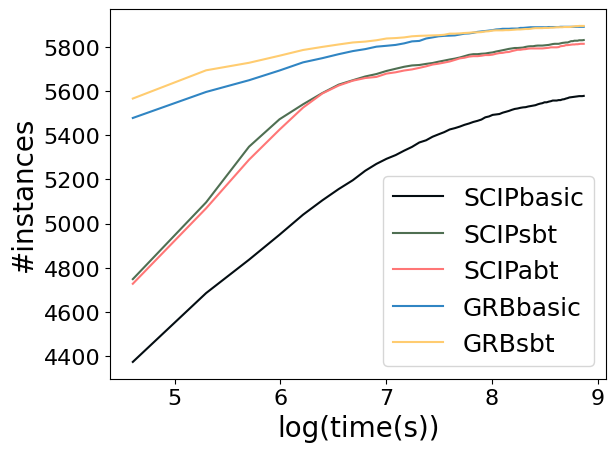}
        \end{subfigure}
        \hfill
        \begin{subfigure}[b]{0.32\textwidth}
            \centering
            \includegraphics[width=\textwidth]{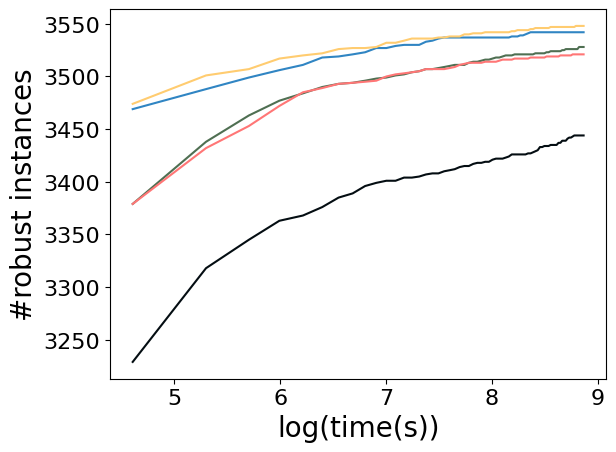}
        \end{subfigure}
        \hfill
        \begin{subfigure}[b]{0.315\textwidth}
            \centering
            \includegraphics[width=\textwidth]{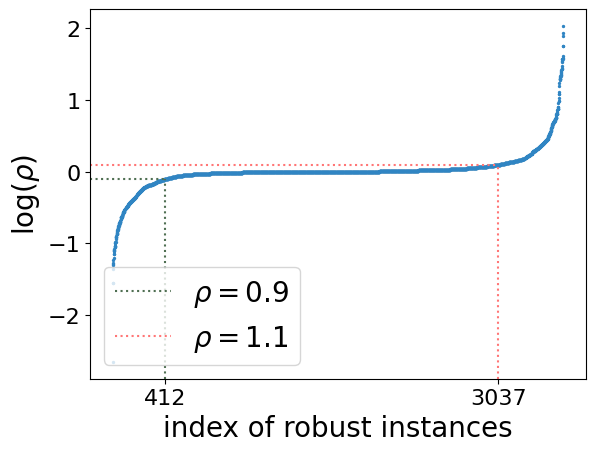}
        \end{subfigure}
        \caption{ENZYMES benchmark. (\textbf{left}) Number of instances solved by each method below different time costs. (\textbf{middle}) Number of robust instances solved by each method below different time costs. (\textbf{right}) Consider $\rho$, the ratio of time cost between SCIPabt and SCIPsbt on each robust instance. SCIPabt is at least $10\%$ faster than SCIPsbt on $412$ robust instances.}
        \label{fig:ENZYMES_time_s2}
    \end{figure*}
    \begin{figure*}[ht]
        \centering
        \begin{subfigure}[b]{0.34\textwidth}
            \centering
            \includegraphics[height=4.3cm]{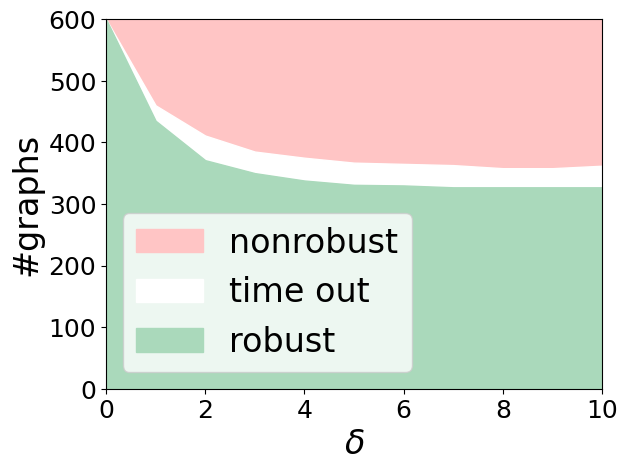}
        \end{subfigure}
        \hfill
        \begin{subfigure}[b]{0.3\textwidth}
            \centering
            \includegraphics[height=4.2cm]{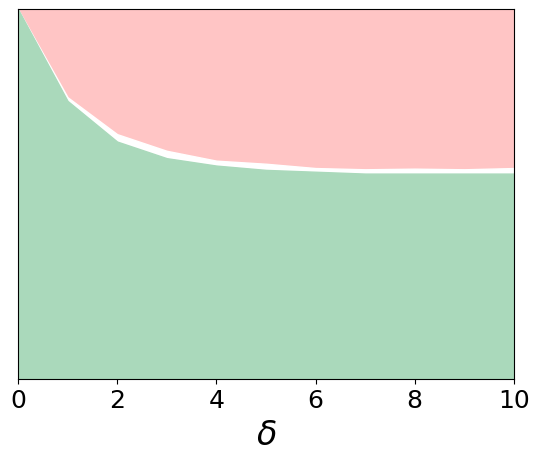}
        \end{subfigure}
        \hfill
        \begin{subfigure}[b]{0.3\textwidth}
            \centering
            \includegraphics[height=4.2cm]{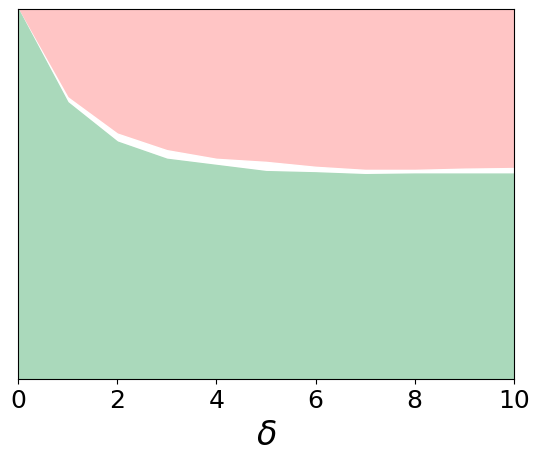}
        \end{subfigure}
        \vskip\baselineskip
        \vspace{-.2in}
        \begin{subfigure}[b]{0.34\textwidth}
            \centering
            \includegraphics[height=4.2cm]{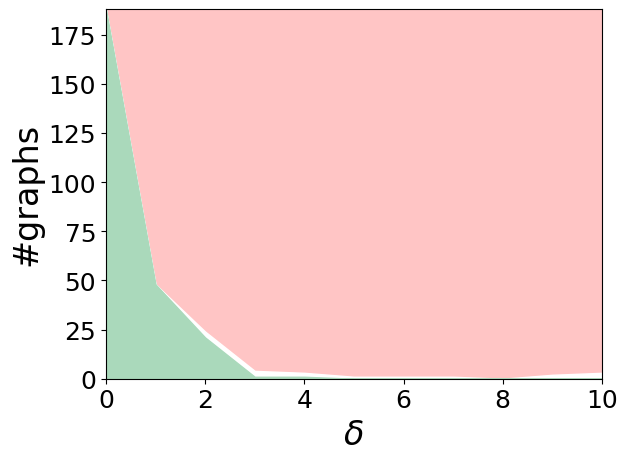}
            \caption{SCIPbasic}
        \end{subfigure}
        \hfill
        \begin{subfigure}[b]{0.3\textwidth}
            \centering
            \includegraphics[height=4.2cm]{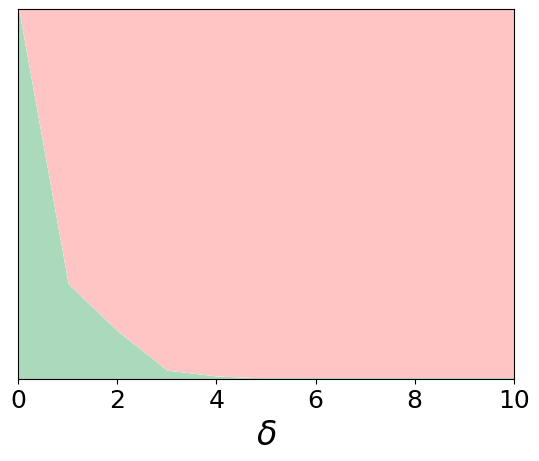}
            \caption{SCIPsbt}
        \end{subfigure}
        \hfill
        \begin{subfigure}[b]{0.3\textwidth}
            \centering
            \includegraphics[height=4.2cm]{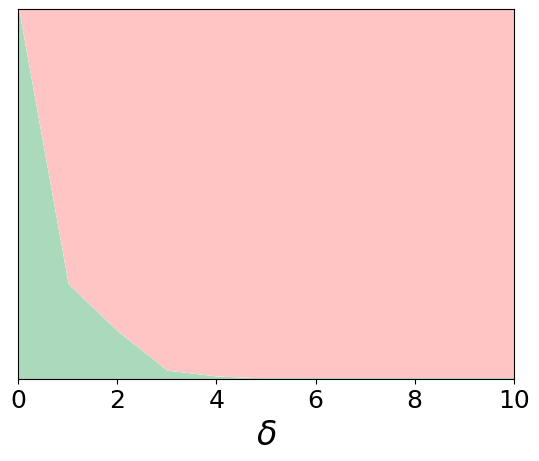}
            \caption{SCIPabt}
        \end{subfigure}
        \vspace{-3mm}
        \caption{For each SCIP-based method with local attack strength $s=2$ on ENZYMES (the first row) and MUTAG (the second row), we count the number of robust graphs (green), nonrobust graphs (red), and time out (white). The percentage $\delta$ of the number of edges is the global budget.}
        \label{fig:visualization_grpah_classification}
    \end{figure*}
    
    This section empirically evaluates the impact of static and aggressive bounds tightening on verifying MPNNs by solving a mixed-integer program (MIP) as shown in Section \ref{sec:GNN_encoding} and Section \ref{sec:verification}.  All GNNs are built and trained using PyG (PyTorch Geometric) 2.1.0 \cite{Fey2019}. All MIPs are implemented in C/C++ using the open-source MIP solver SCIP 8.0.4 \cite{Bestuzheva2023}; all LP relaxations are solved using Soplex 6.0.4 \cite{Gamrath2020}. We used the GNN pull request \cite{Zhang2023b} in the Optimization and Machine Learning Toolkit OMLT \cite{Ceccon2022} to debug the implementation. Observe that we could have alternatively extended a similar tool in SCIP \cite{Turner2023}. For each verification problem, we apply the basic implementation (SCIPbasic), static bounds tightening (SCIPabt), and aggressive bounds tightening (SCIPsbt). Our experiments also include a basic and static bounds tightening implementation of Gurobi 10.0.3 (GRBbasic, GRBsbt) \cite{Gurobi2023}. The code is available at \href{https://github.com/christopherhojny/SCIP-MPNN}{GitHub}, also see \citet{zenodo}.

\subsection{Implementation details}
    Our code models and solves the verification problem in three stages (for basic and static bounds tightening) or four stages (for aggressive bounds tightening). First, parameters of a trained MPNN (weights, biases) and the verification problem (predictive label $c^*$, attack label $c$, global budget $Q$, local budget $q_v$) are read. Second, lower and upper bounds on the variables in Eqs.\ \eqref{eq:big-M of GNN} and \eqref{eq:big-M of ReLU} are computed for SCIPbasic and SCIPsbt. Third, a MIP model of the verification problem is created and solved. We do not solve the MIP models to global optimality as we only ask: \emph{Is this instance robust or not?} We interrupt the optimization once a solution with negative objective value is found, i.e., the instance is non-robust, or the dual bound of the branch-and-bound tree is positive, i.e., the instance is robust.

    In case of SCIPabt, the model is created as for static bounds tightening. During the solving process, a fourth step takes place. This step collects, at each node of the SCIP branch-and-bound tree, the $A$-variables that have been fixed to 0 and 1. We then iterate through the layers of the MPNN and, for each layer, we: (i) recompute the variable bounds as in the second step for the current layer, taking the fixed variables into account following Section \ref{subsec:abt}; (ii) check if any inequality in Eqs.\ \eqref{eq:big-M of GNN} and \eqref{eq:big-M of ReLU} using the newly computed bounds is violated by the current linear programming solution; (iii) if a violated inequality is detected, we add this inequality as a local cutting plane to the model. That is, we inform SCIP that the inequality is only valid at the current node of the branch-and-bound tree and all its children, which is necessary as the inequality is based on fixed variables at the current node. The separation of local cutting planes has been implemented via a separator callback in SCIP.

    To compare with the interval arithmetic approaches, we implemented an OBBT routine. For each variable, we changed the objective to the variable's lower or upper bound, relaxed all binary variables, and solved the resulting LP.

    Besides our implementation in SCIP, we also conducted experiments using Gurobi 10.0.3 \cite{Gurobi2023} to compare baseline (GRBbasic) and static bounds tightening (GRBsbt). We used our SCIP implementation to create the MIP models, store them on the hard drive, and read them via Gurobi's Python interface. We did not investigate aggressive bounds tightening in Gurobi as Gurobi does not support local cutting planes. As for SCIP, we interrupt the solving process after deriving (non-) robustness.

\subsection{Experimental setup}
    All experiments have been conducted on a Linux cluster with 12 Intel Xeon Platinum 8260 2.40 GHz processors each having 48 physical threads. Every model has been solved (either by SCIP or Gurobi) using a single thread. Due to the architecture of the cluster, the jobs have not been run exclusively. The reported running time in the following only consists of the time solving a model, but not creating it. That is, the time for computing the initial variable bounds for baseline and static bounds tightening are not considered, whereas the time for computing bounds in aggressive bounds tightening is considered since this takes place dynamically during the solving process.

    We evaluate the performance of various verification methods on benchmarks including: (i) MUTAG and ENZYMES \cite{Morris2020} for graph classification, and (ii) Cora and CiteSeer \cite{Yang2023} for node classification. All datasets are available in PyG and summarized in Table \ref{tab:benchmarks}. The attack label for each graph/node is fixed as $c = (c^*+1)\bmod C$, where $c^*$ is the predictive label.

    For graph classification, we train a MPNN with $L=3$ SAGEConv \cite{Hamilton2017} layers with $16$ hidden channels, followed by an add pooling and a dense layer as the final classifier. ReLU is used in the first $2$ SAGEConv layers. Following \citet{Jin2020}, $30\%$ of the graphs are used to train the model. The local budget of each node is $q_v=\max\{0,d_v-\max\limits_{u\in V}d_u+s\}$, where $d_v$ is the degree of node $v$, $s$ is the so-called local attack strength. In our experiments, we choose $s\in\{2,3,4\}$, and use $\delta$ percentage of the number of edges as the global budget $Q$, where $1\le \delta\le 10$. We report results for $s=2$ in the paper and results for $s\in\{3,4\}$ in Appendix \ref{app:full_results}.

    For node classification, we train a MPNN with $L=2$ SAGEConv layers with $32$ hidden channels. Similar to \citet{Zugner2020}, $10\%$ of the nodes are used for training. We set $10$ as the global budget and $5$ as the local budget. For each node, we extract its $2$-hop neighborhood to build the corresponding verification problem. It is noteworthy that $2$-hop neighborhood is sufficient for $2$ message passing steps in MPNN, but insufficient for $2$ graph convolutional steps in GCN. The reason is that removing an edge within a $3$-hop neighborhood may influence the degree of a node within a $2$-hop neighborhood.

    All models are trained $200$ epochs with learning rate $0.01$, weight decay $10^{-4}$, and dropout $0.5$. We set $2$ hours as the time limit for verifying a graph in graph classification, and $30$ minutes for verifying a node in node classification.

\subsection{Numerical results}\label{subsec:numerical results}
    For each verification problem, we will get one of three results: robust (objective has non-negative lower bound), nonrobust (a feasible attack, i.e., solution with negative objective, is found), or time out (inconclusive within time limit). For each benchmark, we consider three criteria for each method: (i) average solving time (avg-time), (ii) shifted geometric mean (sgm-time) of solving time, and (iii) number of solved instances within time limit. A commonly used measure to compare MIP-based methods, the shifted geometric mean of $t_1,\cdots,t_n$ is $\left(\prod_{i=1}^n{(t_i+s)}\right)^{1/n}-s$, where shift $s = 10$. Since all approaches use the same model, we ignore model creation time ($\sim 0.036$s). We also exclude the negligible time ($\sim 0.005$s) spent computing variable bounds for both \emph{basic} and \emph{sbt}. For \emph{abt}, we include the time calculating variable bounds since tightening happens at each branch-and-bound tree node. We cannot know the ground truth of a robust instance without complete enumeration or relying on a solver's numerical tolerances, so we classify an instance as robust if all methods claim it is robust except for time out. Three criteria are evaluated for each method on each benchmark over all instances and all robust instances, respectively. We exclude all instances with inconsistencies, i.e., SCIP and Gurobi declares differently, for a fair comparison. See the Appendix \ref{app:implementation_details} for more details.

    The results for node classification are reported in Table \ref{tab:summary_node_classification} in Appendix \ref{app:full_results}. Only removing edges results in simple verification problems: all methods can solve all instances instantly. Adding more cutting planes is not helpful as this could hinder heuristics to find a feasible attack (in case of non-robustness) or result in solving more (difficult) LPs due to additional cutting planes (in case of robustness). Therefore, we skip the comparison with GRBbasic and GRBsbt. 
    
    For graph classification, we analyze results  with local attack strength $s=2$ in the main text and put results for $s\in\{3,4\}$ in Appendix \ref{app:full_results}. Our numerical analysis for $s=2$ is consistent with the $s\in\{3,4\}$ results. Table \ref{tab:summary_graph_classification_s2} summarizes the results for all instances. Figure \ref{fig:ENZYMES_time_s2} visualizes the number of solved (robust) instances below different time costs, and compares the time costs between SCIPabt and SCIPsbt for ENZYMES. See Appendix \ref{app:full_results} for similar plots for MUTAG. Figure \ref{fig:visualization_grpah_classification} plots the number of different types of graphs (robust, nonrobust, time out) with various global budgets $1\le \delta\le 10$ for ENZYMES and MUTAG. 

    As shown in Table \ref{tab:summary_graph_classification_s2}, SCIPsbt is around three times faster than SCIPbasic and solves more instances within the same time limit. From the comparison between GRBbasic and GRBsbt, we can still notice the speed-up from static bounds tightening. Considering all instances from MUTAG, it seems like static bounds tightening slows down the solving process. The reason is that most MUTAG instances are not robust, i.e., finding good bounds on the objective value is unnecessary, finding a feasible attack instead is sufficient. 

    Our numerical results reflect similar performance between SCIPsbt and SCIPabt w.r.t.\ times in Table \ref{tab:summary_graph_classification_s2}. SCIPabt might even be slower than SCIPsbt in some instances. On the one hand, we proposed \emph{abt} as a general extension of \emph{sbt} and expect it outperforms \emph{sbt} in harder verification problems. One can easily create larger problems from different aspects, e.g., increasing budgets, incorporating feature perturbations, and enlarging size of models. However, larger problems do not imply harder verification problems since the instance could be nonrobust. Then the extra cutting planes added in \emph{abt} may slow down finding a feasible attack. On the other hand, as shown in Appendix \ref{app:full_results}, when comparing OBBT and SCIPsbt over robust instances, OBBT bounds are indeed tighter but still perform similarly to SCIPsbt. The phenomenon that tighter bounds can result in slower solving times has been previously reported \citep{Badilla2023}.

    Observe in Figure \ref{fig:ENZYMES_time_s2} that, of the 3549 robust instances from the ENZYMES benchmark, 412 are significantly faster when using SCIPabt and 512 are significantly faster when using SCIPsbt. Similarly for the MUTAG benchmark, 19 of the 44 robust instances are substantially faster using SCIPabt rather than SCIPsbt (see Figure \ref{fig:MUTAG_time_full}). This is also reflected in Table \ref{tab:summary_graph_classification_s2} when considering the robust instances only. For the MUTAG instances (which are harder to solve than ENZYMES), SCIPabt is roughly 10\% faster than SCIPsbt in shifted geometric mean (29\% in arithmetic mean). This effect is even more pronounced for the most difficult MUTAG instances with a global budget of 2\% and 3\%, where SCIPabt is 10.8\% and 22.0\%, respectively, faster than SCIPsbt in shifted geometric mean (29.4\% and 35.8\% in arithmetic mean), see Table \ref{tab:MUTAG_s2_full} in Appendix \ref{app:full_results}.

    We therefore propose running SCIPsbt and SCIPabt in parallel, this idea corresponds to the common observation in MIP that parallelizing multiple strategies (here: SCIPsbt and SCIPabt) yields more benefits than parallelizing just one algorithm \cite{Carvajal2014}. 

\section{Conclusion}\label{sec:conclusion}
    We propose topology-based bounds tightening approaches to verify message-passing neural networks. By exploiting graph structures and available budgets, our techniques compute tighter bounds for variables and thereby help certify robustness. Numerical results show the improvement of topology-based bounds tightening w.r.t.\ the solving time and the number of solved instances.  

\section*{Acknowledgements}
    This work was supported by the Engineering and Physical Sciences Research Council [grant number EP/W003317/1], an Imperial College Hans Rausing PhD Scholarship to SZ, and a BASF/RAEng Research Chair in Data-Driven Optimisation to RM.
    This collaboration was initiated at the Mittag-Leffler Institute seminar titled \emph{Learning from Both Sides: Linear and Nonlinear Mixed-Integer Optimization}.

\section*{Impact Statement}
    Despite the widespread use and outstanding performance of GNNs in various fields, their vulnerabilities are frequently detected, even with slight perturbations. In safety-critical applications such as drug design and autonomous driving, such vulnerabilities could result in unacceptable consequences. For safety considerations, two major directions are investigated by many researchers: (i) verification: how to measure the robustness of GNNs, and (ii) robust training: how to train GNNs that are more robust. This work considers verification problems on MPNNs, a classic GNN framework whose robustness is seldom studied in the literature. With the proposed bounds tightening strategies, we hope that the robustness of MPNNs deployed in real-world applications could be verified efficiently to avoid unreliable predictions.


\bibliography{paper.bib}
\bibliographystyle{icml2024}

\newpage
\appendix
\onecolumn

\section{Proofs of Properties \ref{prop:sbt} \& \ref{prop:abt}}\label{app:proofs}
    \begin{proof}[\textbf{Proof of Property \ref{prop:sbt}}]
        Since $A_{u,v}\ge 0$, we only need to consider the lower bound of $\sum\limits_{f\in F_{l-1}}w_{f,f'}x_{u,f}$ in Eq. \eqref{eq:simplified definition of lower bound}, which is denoted as $lb_{u\to v}$. Using bounds for $x_{u,f}$ gives $lb_{u\to v}$ defined as Eq. \eqref{eq:contribution of neighbors}. Then Eq. \eqref{eq:simplified definition of lower bound} becomes:
        \begin{equation}
            \begin{aligned}
                lb(\bar x_{v,f'})=\min\limits_{A\in\mathcal P_1(A^*)}\sum\limits_{u\in V}A_{u,v}lb_{u\to v}+b_{v,f'}.
            \end{aligned}
        \end{equation}
        Recall the definition of $\mathcal P_1(A^*)$, we can remove/add at most $q_v$ neighbors of node $v$. Denote the set of removed/added neighbors as $V_{lb}$. Then the summation $u\in V$ can be divided into: $u\in \mathcal N^*(v)$ (original neighbors), $u\in V_{lb}\backslash \mathcal N^*(v)$ (added neighbors), $u\in \mathcal N^*(v)\cap V_{lb}$ (removed neighbors), and $u\in V\backslash (\mathcal N^*(v)\cup V_{lb})$ (nodes without contribution either before or now), from which we obtain that:
        \begin{equation}
            \begin{aligned}
                \sum\limits_{u\in V}A_{u,v}lb_{u\to v}&=\left(\sum\limits_{u\in \mathcal N*(v)}+\sum\limits_{u\in V_{lb}\backslash\mathcal N*(v)}-\sum\limits_{u\in \mathcal N^*(v)\cap V_{lb}}\right)lb_{u\to v}\\
                &=\sum\limits_{u\in\mathcal N^*(v)}lb_{u\to v} + \sum\limits_{u\in V_{lb}}\Delta_{u\to v}
            \end{aligned}
        \end{equation}
        with $\Delta_{u\to v}$ defined as Eq. \eqref{eq:changes of add/remove neighbors}. Moving the minimization to the only changeable term $\sum\limits_{u\in V_{lb}}\Delta_{u\to v}$ yields Eq. \eqref{eq:decomposed definition of lower bound}.
    \end{proof}
    
    \begin{proof}[\textbf{Proof of Property \ref{prop:abt}}]
        This property can be proved similarly to Property \ref{prop:sbt}. Here we give a simple way to prove it using the conclusion of Property \ref{prop:sbt}. As mentioned in Remark \ref{remark:abt is sbt on each node}, we can modify the original graph with previous decisions, i.e., current neighborhood of node $v$ is $(\mathcal N^*(v)\backslash V_0)\cup V_1$ and the remaining budget is $q'_v$ defined as Eq. \eqref{eq:budgets in B&B}. Since we already decided the values of $A_{u,v}$ for $u\in V_0\cup V_1$, $V_{lb}$ should not include any node in $V_0\cup V_1$. Applying Property \ref{prop:sbt} on the current graph with budgets $q'_v$ gives Eq. \eqref{eq:decomposed lower bound in B&B} and finishes this proof. 
    \end{proof}

\section{Implementation details}\label{app:implementation_details}
\subsection{Local cutting planes}
    A local cutting plane in SCIP is any inequality together with the first node where it is valid. As shown in the right side of Figure \ref{fig:overall}, after branching variable $A_{1,0}=0$, we can improve the bounds of node $0$ from $[-4,9]$ to $[-3,6]$. These improved bounds can be incorporated into the model by adding inequalities $x\ge -3$ and $x\le 6$ to the left branch. However, these inequalities cannot be added to the root node as they are invalid unless $A_{1,0}=0$. At each node in the branch-and-bound tree, all local cutting planes associated with its ancestors are valid. Including all of these local cutting planes is correct but inefficient. Therefore, SCIP will decide whether an inequality will be added. For example, if node $0$ equals to $7$ at the left child of the root node, then SCIP only adds $x\le 6$ into the model since this inequality is violated.

\subsection{Inconsistent instances}
    We exclude instances with inconsistent results from different solvers, which results in the numbers in column ``\#" of Table \ref{tab:ENZYMES_s2_full} -- Table \ref{tab:MUTAG_s4_full} being smaller than the total number of graphs from each benchmark. This is because we found and reported a small bug in Gurobi wherein several instances were declared \emph{infeasible} despite having a feasible solution found by SCIP. This bug was easily fixed in Gurobi and is now incorporated into a recent minor release. We continued using the Gurobi version with the bug after the ICML rebuttal period because there was insufficient time to redo all the experiments. So, in fairness to the Gurobi solver, we slightly modified the optimization problems by increasing the right-hand side of Eq.~\eqref{eq:big-M of ReLU3} by $10^{-11}$. This change eliminated the bug, but then meant that the Gurobi solver frequently returns solutions that are infeasible, i.e., violating some constraints of the MIP model by more than tolerance $10^{-6}$. To compare on an equal footing, we only report results on instances where all solvers return consistent answers. In other words, if one solver wrongly declares an instance \emph{infeasible} or returns an infeasible solution, we do not report results.

\newpage
\section{Full numerical results}\label{app:full_results}
\begin{figure*}[h]
    \centering
    \begin{subfigure}[c]{0.34\textwidth}\raggedleft
        \myrowlabel{$s=2$}
        \raisebox{-.5\height}{\includegraphics[height=4cm]{Figures/ENZYMES_SCIPbasic_s2.png}}\\
        \myrowlabel{$s=3$}
        \raisebox{-.5\height}{\includegraphics[height=4cm]{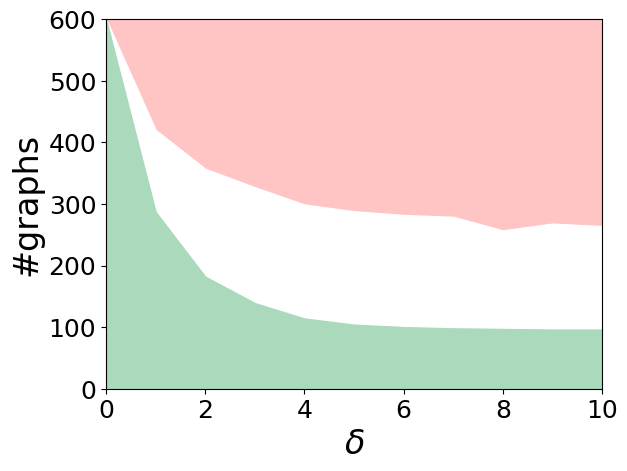}}\\
        \myrowlabel{$s=4$}
        \raisebox{-.5\height}{\includegraphics[height=4cm]{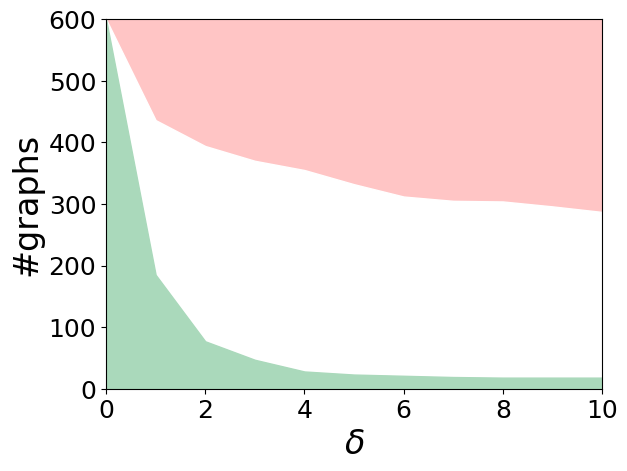}}
    \caption{SCIPbasic}
    \end{subfigure}%
    \hspace{1em}
    \begin{subfigure}[c]{0.30\textwidth}\raggedleft
        \includegraphics[height=4cm]{Figures/ENZYMES_SCIPsbt_s2.png}
        \includegraphics[height=4cm]{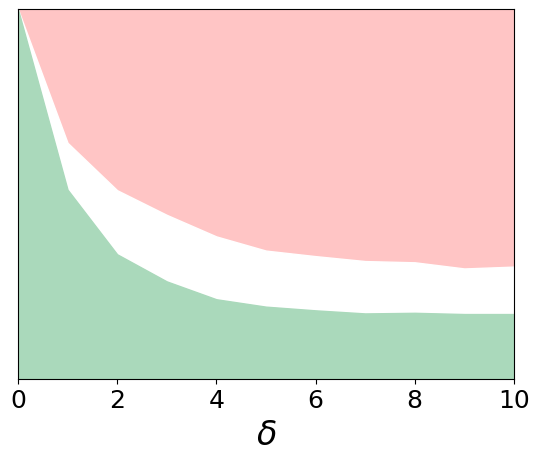}
        \includegraphics[height=4cm]{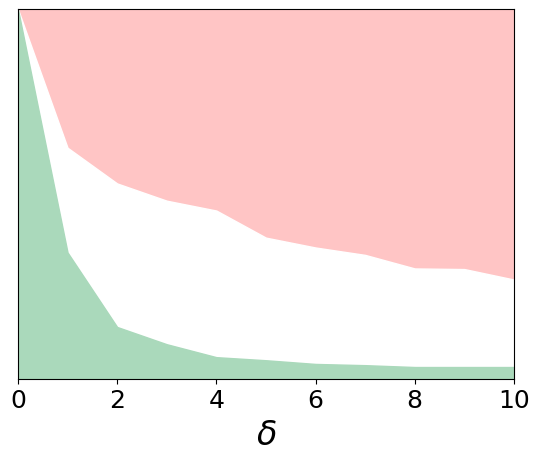}
    \caption{SCIPsbt}
    \end{subfigure}
    \hspace{1em}
    \begin{subfigure}[c]{0.30\textwidth}\centering
        \includegraphics[height=4cm]{Figures/ENZYMES_SCIPabt_s2.png}
        \includegraphics[height=4cm]{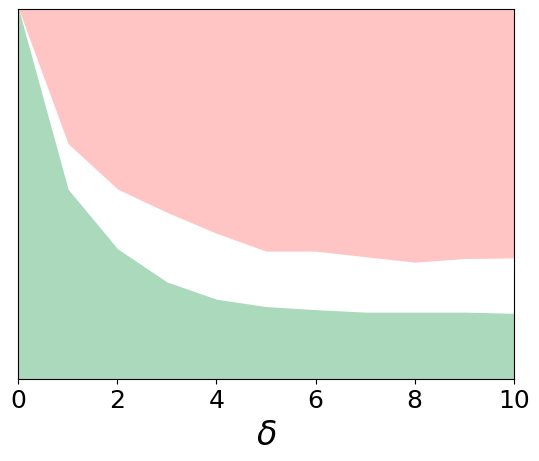}
        \includegraphics[height=4cm]{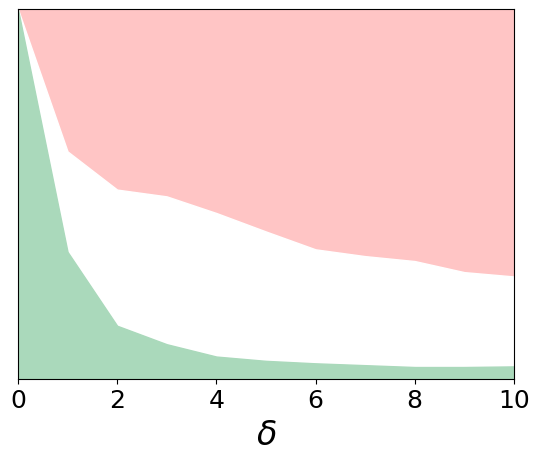}
    \caption{SCIPabt}
    \end{subfigure}
    \caption{For each SCIP-based method on ENZYMES with local attack strength $s\in \{2,3,4\}$, we count the number of robust graphs (green), nonrobust graphs (red), and time out (white). The percentage $\delta$ of the number of edges is the global budget.}
    \label{fig:ENZYMES_budget_full}
\end{figure*}

\begin{figure*}[h]
    \centering
    \begin{subfigure}[c]{0.34\textwidth}\raggedleft
        \myrowlabel{$s=2$}
        \raisebox{-.5\height}{\includegraphics[height=4cm]{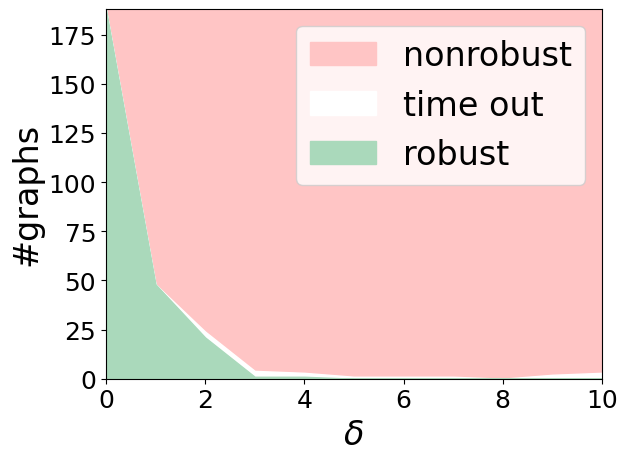}}\\
        \myrowlabel{$s=3$}
        \raisebox{-.5\height}{\includegraphics[height=4cm]{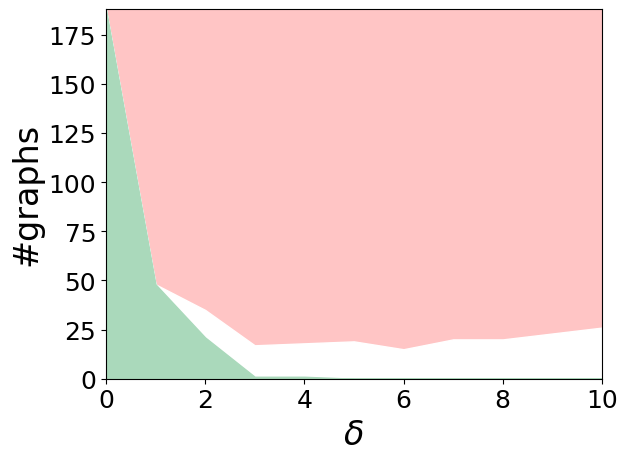}}\\
        \myrowlabel{$s=4$}
        \raisebox{-.5\height}{\includegraphics[height=4cm]{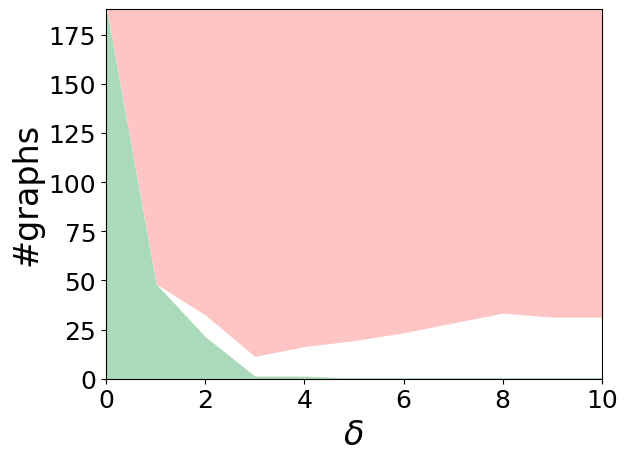}}
        \caption{SCIPbasic}
    \end{subfigure}%
    \hspace{1em}
    \begin{subfigure}[c]{0.30\textwidth}\raggedleft
        \includegraphics[height=4cm]{Figures/MUTAG_SCIPsbt_s2.png}
        \includegraphics[height=4cm]{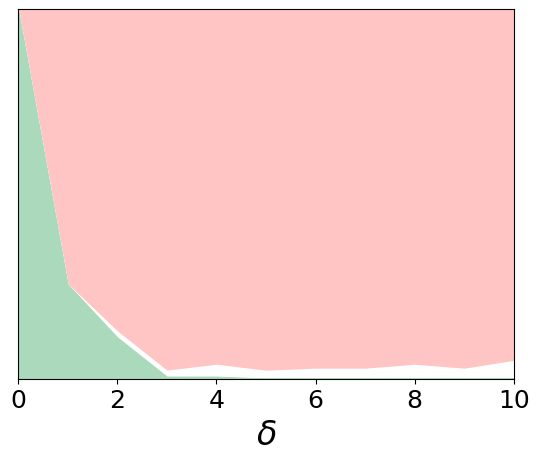}
        \includegraphics[height=4cm]{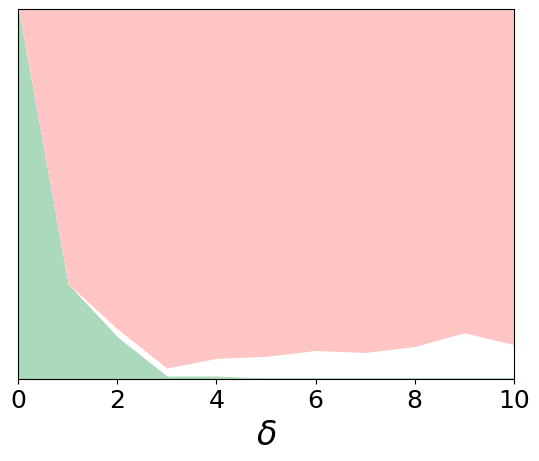}
    \caption{SCIPsbt}
    \end{subfigure}
    \hspace{1em}
    \begin{subfigure}[c]{0.30\textwidth}\centering
        \includegraphics[height=4cm]{Figures/MUTAG_SCIPabt_s2.png}
        \includegraphics[height=4cm]{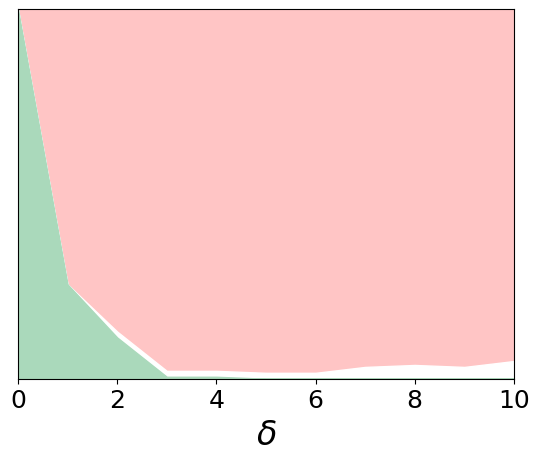}
        \includegraphics[height=4cm]{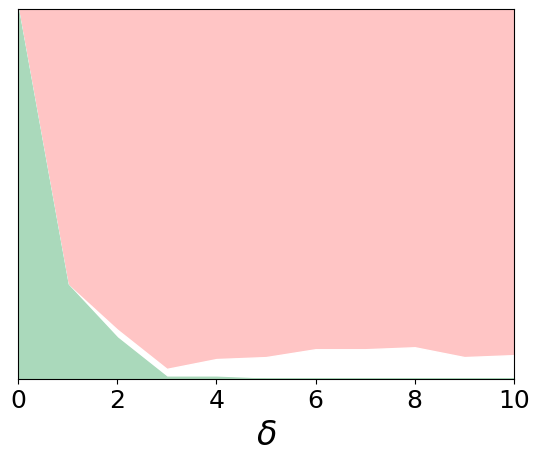}
    \caption{SCIPabt}
    \end{subfigure}
    \caption{For each SCIP-based method on MUTAG with local attack strength $s\in \{2,3,4\}$, we count the number of robust graphs (green), nonrobust graphs (red), and time out (white). The percentage $\delta$ of the number of edges is the global budget.}
    \label{fig:MUTAG_budget_full}
\end{figure*}

\begin{figure*}[h]
    \centering
    \begin{subfigure}[c]{0.34\textwidth}\raggedleft
        \myrowlabel{$s=2$}
        \raisebox{-.5\height}{\includegraphics[height=4cm]{Figures/ENZYMES_all_instances_s2.png}}\\
        \myrowlabel{$s=3$}
        \raisebox{-.5\height}{\includegraphics[height=4cm]{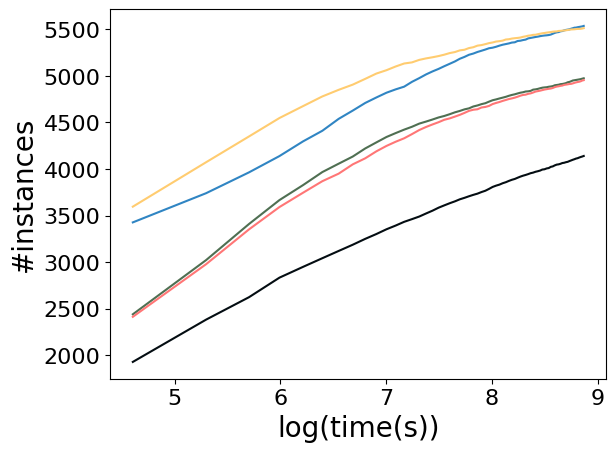}}\\
        \myrowlabel{$s=4$}
        \raisebox{-.5\height}{\includegraphics[height=4cm]{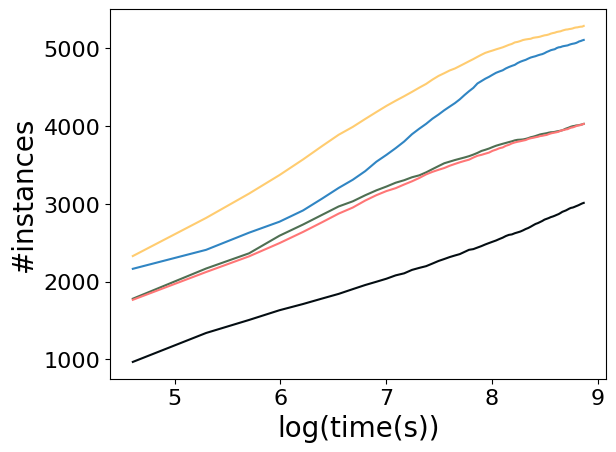}}
    \end{subfigure}%
    \hspace{1em}
    \begin{subfigure}[c]{0.30\textwidth}\raggedleft
        \includegraphics[height=4cm]{Figures/ENZYMES_robust_instances_s2.png}
        \includegraphics[height=4cm]{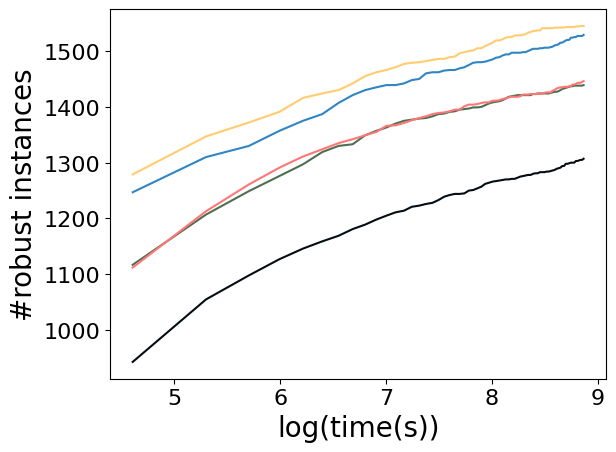}
        \includegraphics[height=4cm]{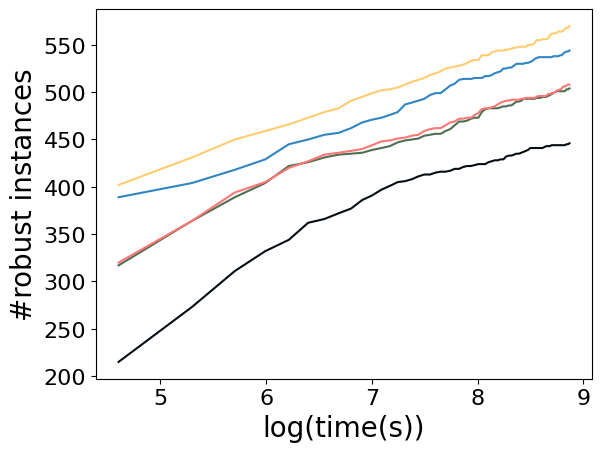}
    \end{subfigure}
    \hspace{1em}
    \begin{subfigure}[c]{0.30\textwidth}\centering
        \includegraphics[height=4cm]{Figures/ENZYMES_ratio_s2.png}
        \includegraphics[height=4cm]{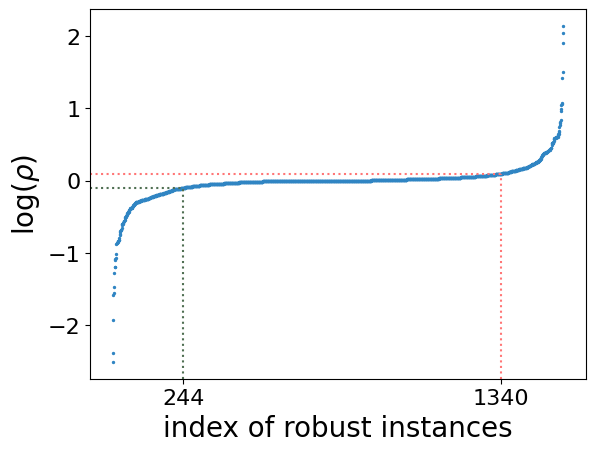}
        \includegraphics[height=4cm]{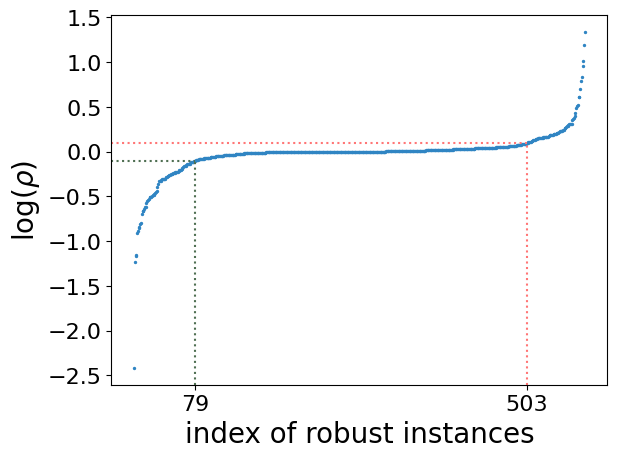}
    \end{subfigure}
    \caption{ENZYMES benchmark with local attack strength $s\in\{2,3,4\}$. (\textbf{left}) Number of instances solved by each method below different time costs. (\textbf{middle}) Number of robust instances solved by each method below different time costs. (\textbf{right}) Consider $\rho$, the ratio of time cost between SCIPabt and SCIPsbt on each robust instance.}
    \label{fig:ENZYMES_time_full}
\end{figure*}

\begin{figure*}[h]
    \centering
    \begin{subfigure}[c]{0.34\textwidth}\raggedleft
        \myrowlabel{$s=2$}
        \raisebox{-.5\height}{\includegraphics[height=4cm]{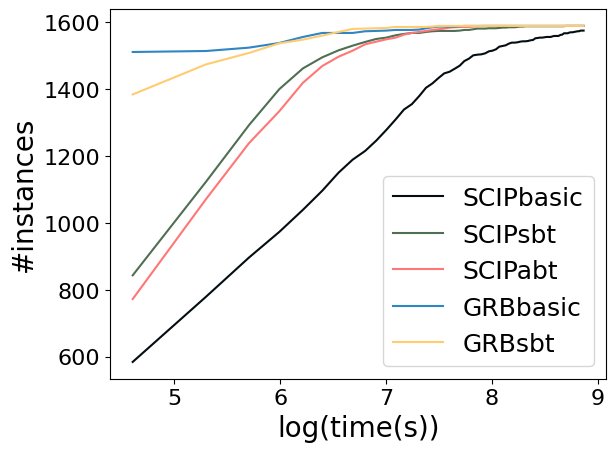}}\\
        \myrowlabel{$s=3$}
        \raisebox{-.5\height}{\includegraphics[height=4cm]{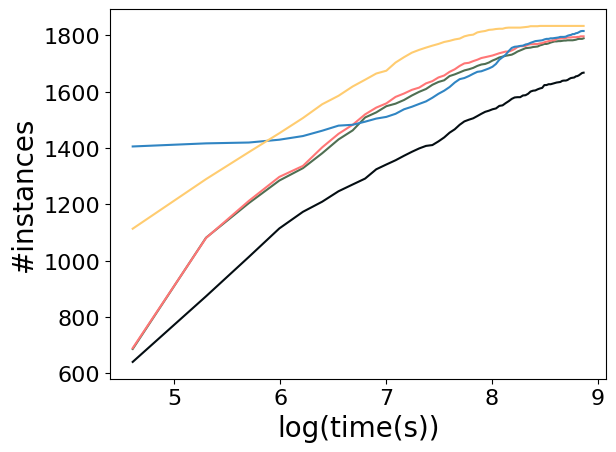}}\\
        \myrowlabel{$s=4$}
        \raisebox{-.5\height}{\includegraphics[height=4cm]{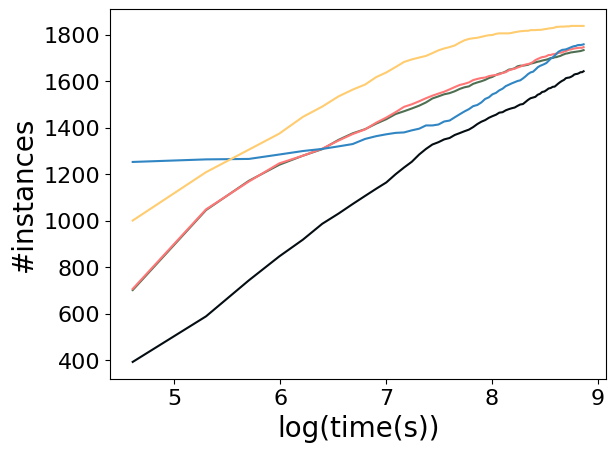}}
    \end{subfigure}%
    \hfill
    \begin{subfigure}[c]{0.30\textwidth}\raggedleft
        \includegraphics[height=4cm]{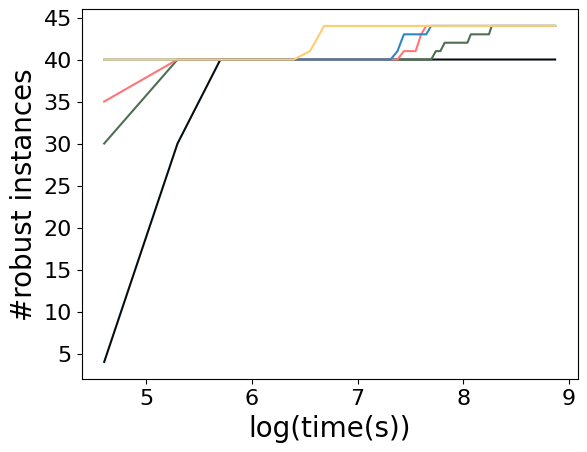}
        \includegraphics[height=4cm]{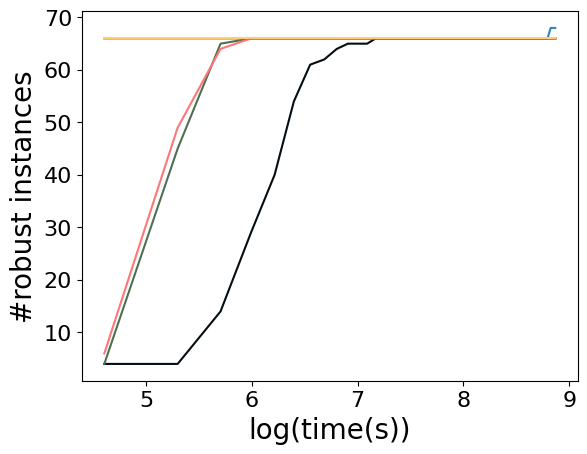}
        \includegraphics[height=4cm]{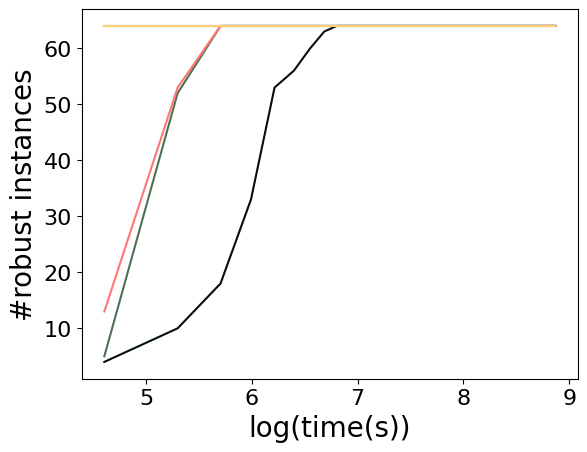}
    \end{subfigure}
    \hfill
    \begin{subfigure}[c]{0.30\textwidth}\centering
        \includegraphics[height=4cm]{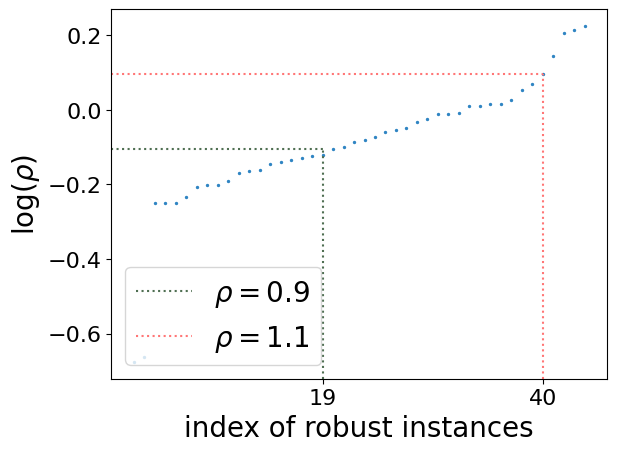}
        \includegraphics[height=4cm]{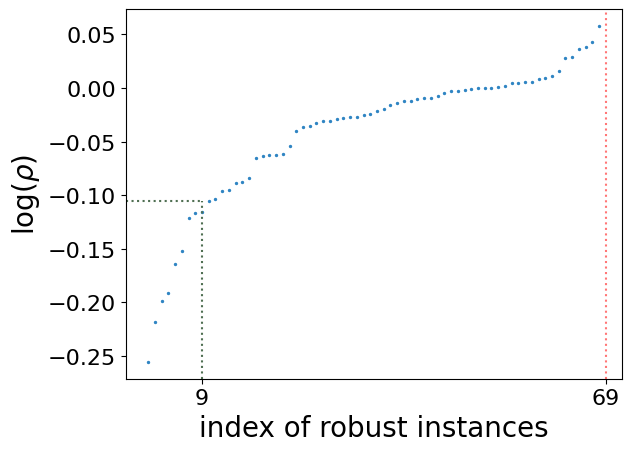}
        \includegraphics[height=4cm]{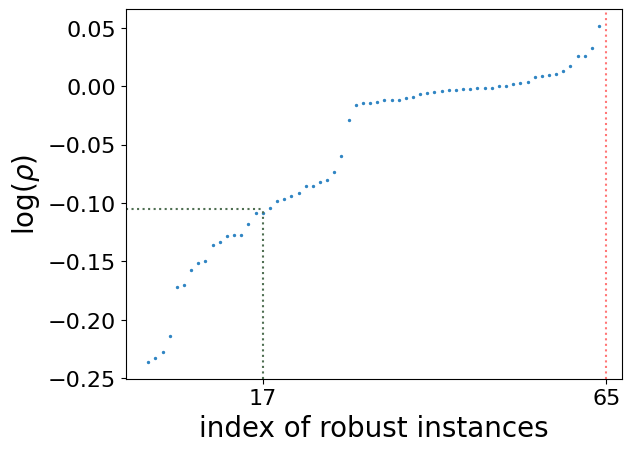}
    \end{subfigure}
    \caption{MUTAG benchmark with local attack strength $s\in\{2,3,4\}$. (\textbf{left}) Number of instances solved by each method below different time costs. (\textbf{middle}) Number of robust instances solved by each method below different time costs. (\textbf{right}) Consider $\rho$, the ratio of time cost between SCIPabt and SCIPsbt on each robust instance.}
    \label{fig:MUTAG_time_full}
\end{figure*}

    Table \ref{tab:summary_node_classification} gives results for node classification. Tables \ref{tab:ENZYMES_full} and \ref{tab:MUTAG_full} summarize the results for graph classification with local attack strength $s\in\{2,3,4\}$ on ENZYMES and MUTAG, respectively. Detailed results with different global budgets are reported in Tables \ref{tab:ENZYMES_s2_full}--\ref{tab:MUTAG_s4_full}. Figures \ref{fig:ENZYMES_time_full} and \ref{fig:MUTAG_time_full} plot the number of different types of graphs (robust, nonrobust, time out) with various budgets. Figures \ref{fig:ENZYMES_budget_full} and  \ref{fig:MUTAG_budget_full} visualize the number of solved (robust) instances below different time costs, and compares the time costs between SCIPabt and SCIPsbt. For graph classification, our numerical analysis in Section \ref{subsec:numerical results} for $s=2$ is consistent with results for $s\in\{3,4\}$, as shown here.

    Furthermore, we report some numerical results to show that applying OBBT in GNN verification is not recommended. Our OBBT implementation is straight-forward: for each variable, change the objective in our verification problem to its lower/upper bound and relax all binary variables into continuous variables. Due to the high time cost for solving each instance, as well as our numerical observation that tighter bounds may slow down the process to find a feasible attack for a non-robust instance, we only apply OBBT to  the 44 robust instances for MUTAG with $s=2$. The average time spent on calculating OBBT bounds for all variables is 6045.02s (3.57s per LP), and the average solving time with OBBT bounds is 331.60s. As shown in Table \ref{tab:summary_graph_classification_s2}, the average solving time for SCIPsbt is 336.41s. Therefore, the improvement of OBBT w.r.t.\ the solving time is quite limited, despite its high time cost on computing bounds. We also compare the bounds derived from \emph{sbt} and OBBT using the relative bound tightness (RBT) used in \citet{Badilla2023}: $RBT=\frac{SBT-OBBT}{OBBT+10^{-10}}$, where SBT/OBBT represents the length of interval bounds derived from \emph{sbt}/OBBT, respectively. For each MPNN layer, we average RBT over all variables in this layer for all 44 instances. The resulting RBT values are $0.07,0.22,0.57$, which means that the \emph{sbt} bounds are quite closed to OBBT bounds for early layers and decently tight for deeper layers.

\begin{table*}[h]
    \caption{Summary of results for node classification. 
    The approaches tested are SCIP basic (SCIPbasic), SCIP static bounds tightening (SCIPsbt), and SCIP aggressive bounds tightening (SCIPabt). For each global budget, the number of instances is 2708 for Cora and 3312 for CiteSeer. We compare times with respect to both average (``avg-time") and the shifted geometric mean (``sgm-time"), as well as the number of solved instances within time limits 30 minutes (``\# solved"). Since robust instances are the ones where mixed-integer performance is most important (non-robust instances may frequently be found by more heuristic approaches), we also compare on the set of robust instances.
    }
    \label{tab:summary_node_classification}
    \vskip 0.15in
    \begin{center}
    \begin{small}
    \begin{tabular*}{\textwidth}{@{\;\;\extracolsep{\fill}}llrrrrrrrr}
        \toprule
        \multirow{2}{*}{benchmark} & \multirow{2}{*}{method} & \multicolumn{4}{c}{all instances} & \multicolumn{4}{c}{robust instances}\\
        \cmidrule{3-6} \cmidrule{7-10}
         & & \# & avg-time(s) & sgm-time(s) & \# solved  & \# & avg-time(s) & sgm-time(s) & \# solved\\
         \midrule
        \multirow{3}{*}{Cora} 
        &SCIPbasic     &  2708 &     \textbf{0.10} &     \textbf{0.10} &  \textbf{2708} &  2223 &     \textbf{0.08} &     \textbf{0.08} &  \textbf{2223} \\
        &SCIPsbt       &  2708 &     0.17 &     0.16 &  \textbf{2708} &  2223 &     0.10 &     0.10 &  \textbf{2223} \\
        &SCIPabt       &  2708 &     0.46 &     0.38 &  \textbf{2708} &  2223 &     0.16 &     0.14 &  \textbf{2223} \\
        \midrule
        \multirow{3}{*}{CiteSeer} 
        & SCIPbasic     &  3312 &     \textbf{0.07} &     \textbf{0.06} &  \textbf{3312} &  2917 &     \textbf{0.06} &     \textbf{0.06} &  \textbf{2917} \\
        & SCIPsbt       &  3312 &     \textbf{0.07} &     0.07 &  \textbf{3312} &  2917 &     \textbf{0.06} &     \textbf{0.06} &  \textbf{2917} \\
        & SCIPabt       &  3312 &     0.31 &     0.17 &  \textbf{3312} &  2917 &     0.12 &     0.10 &  \textbf{2917} \\
        \bottomrule
    \end{tabular*}
    \end{small}
    \end{center}
    \vskip -0.1in
\end{table*}

\begin{table*}[h]
    \caption{Summary of results for graph classification on ENZYMES with local attack strength $s\in\{2,3,4\}$. The approaches tested are SCIP basic (SCIPbasic), SCIP static bounds tightening (SCIPsbt), SCIP aggressive bounds tightening (SCIPabt), Gurobi basic (GRBbasic), and Gurobi static bounds tightening (GRBsbt). For each global budget, the number of instances is 600, but we only present comparisons when all methods give consistent results except for time out, as shown in column ``\#". We compare times with respect to both average (``avg-time") and the shifted geometric mean (``sgm-time"), as well as the number of solved instances within time limits 2 hours (``\# solved"). Since robust instances are the ones where mixed-integer performance is most important (non-robust instances may frequently be found by more heuristic approaches), we also compare on the set of robust instances.}
    \label{tab:ENZYMES_full}
    \vskip 0.15in
    \begin{center}
    \begin{small}
    \begin{tabular*}{\textwidth}{@{\;\;\extracolsep{\fill}}llrrrrrrrr}
        \toprule
        \multirow{2}{*}{} & \multirow{2}{*}{method} & \multicolumn{4}{c}{all instances} & \multicolumn{4}{c}{robust instances}\\
        \cmidrule{3-6} \cmidrule{7-10}
         & & \# & avg-time(s) & sgm-time(s) & \# solved  & \# & avg-time(s) & sgm-time(s) & \# solved\\
        \midrule
         \multirow{5}{*}{$s=2$}
            & SCIPbasic  & 5915 & 605.97 & 37.81 & 5579 & 3549 & 278.58 & 12.51 & 3444 \\
            & SCIPsbt    & 5915 & \textbf{230.59} & \textbf{21.26} & \textbf{5831} & 3549 &  \textbf{82.89} & 
    \textbf{6.65} & \textbf{3528} \\
            & SCIPabt    & 5915 & 246.02 & 21.77 & 5817 & 3549 &  88.95 &  6.71 & 3522 \\
            \cmidrule{2-10}
            & GRBbasic   & 5915 &  86.03 &  7.59 & 5892 & 3549 &  32.80 &  2.82 & 3542 \\
            & GRBsbt     & 5915 &  \textbf{74.87} &  \textbf{7.09} & \textbf{5895} & 3549 &  \textbf{22.90} &  \textbf{2.50} & \textbf{3548} \\
        \midrule
        \multirow{5}{*}{$s=3$}  
        & SCIPbasic  & 5855 & 2628.71 &  460.25 & 4149 & 1554 & 1423.02 &   99.43 & 1308 \\
        & SCIPsbt    & 5855 & \textbf{1533.12} &  \textbf{214.04} & \textbf{4977} & 1554 &  750.70 &   48.68 & 1440 \\
        & SCIPabt    & 5855 & 1583.67 &  224.59 & 4943 & 1554 &  \textbf{738.67} &   \textbf{48.11} & \textbf{1447} \\
        \cmidrule{2-10}
        & GRBbasic   & 5855 &  845.77 &   85.82 & \textbf{5549} & 1554 &  367.89 &   22.54 & 1530 \\
        & GRBsbt     & 5855 &  \textbf{738.26} &   \textbf{74.83} & 5524 & 1554 &  \textbf{238.48} &   \textbf{18.89} & \textbf{1545} \\
        \midrule
         \multirow{5}{*}{$s=4$}  
        & SCIPbasic  & 5901 & 4163.96 & 1492.61 & 3003 &  577 & 2030.34 &  319.19 &  448 \\
        & SCIPsbt    & 5901 & \textbf{2802.04} &  \textbf{581.20} & 4012 &  577 & 1353.31 &  154.83 &  504 \\
        & SCIPabt    & 5901 & 2841.31 &  602.57 & \textbf{4013} &  577 & \textbf{1318.26} &  \textbf{152.16} &  \textbf{509} \\
        \cmidrule{2-10}
        & GRBbasic   & 5901 & 1780.68 &  306.11 & 5125 &  577 &  862.25 &   62.75 &  547 \\
        & GRBsbt     & 5901 & \textbf{1372.88} &  \textbf{233.52} & \textbf{5290} &  577 &  \textbf{589.89} &   \textbf{49.66} &  \textbf{570} \\
        \bottomrule
    \end{tabular*}
    \end{small}
    \end{center}
    \vskip -0.1in
\end{table*}

\begin{table*}[h]
    \caption{Summary of results for graph classification on MUTAG with local attack strength $s\in\{2,3,4\}$. The approaches tested are SCIP basic (SCIPbasic), SCIP static bounds tightening (SCIPsbt), SCIP aggressive bounds tightening (SCIPabt), Gurobi basic (GRBbasic), and Gurobi static bounds tightening (GRBsbt). For each global budget, the number of instances is 188, but we only present comparisons when all methods give consistent results except for time out, as shown in column ``\#". We compare times with respect to both average (``avg-time") and the shifted geometric mean (``sgm-time"), as well as the number of solved instances within time limits 2 hours (``\# solved"). Since robust instances are the ones where mixed-integer performance is most important (non-robust instances may frequently be found by more heuristic approaches), we also compare on the set of robust instances.}
    \label{tab:MUTAG_full}
    \vskip 0.15in
    \begin{center}
    \begin{small}
    \begin{tabular*}{\textwidth}{@{\;\;\extracolsep{\fill}}llrrrrrrrr}
        \toprule
        \multirow{2}{*}{} & \multirow{2}{*}{method} & \multicolumn{4}{c}{all instances} & \multicolumn{4}{c}{robust instances}\\
        \cmidrule{3-6} \cmidrule{7-10}
         & & \# & avg-time(s) & sgm-time(s) & \# solved  & \# & avg-time(s) & sgm-time(s) & \# solved\\
        \midrule
       \multirow{5}{*}{$s=2$}
            & SCIPbasic  & 1589 & 679.86 & 189.75 & 1575 & 44 & 798.47 & 202.93 & 40 \\
            & SCIPsbt    & 1589 & \textbf{196.07} &  \textbf{75.17} & \textbf{1589} & 44 & 336.41 & 100.86 & \textbf{44} \\
            & SCIPabt    & 1589 & 207.50 &  82.43 & \textbf{1589} & 44 & \textbf{238.10} &  \textbf{91.11} & \textbf{44} \\
            \cmidrule{2-10}
            & GRBbasic   & 1589 &  \textbf{34.58} &  \textbf{4.06} & \textbf{1589} & 44 & 162.25 &  \textbf{12.11} & \textbf{44} \\
            & GRBsbt     & 1589 &  59.93 & 22.40 & 1589 & 44 &  \textbf{73.78} &  15.00 &  \textbf{44} \\
        \midrule
        \multirow{5}{*}{$s=3$}
        & SCIPbasic  & 1840 & 1355.17 &  262.63 & 1670 &   68 &  643.89 &  427.71 &   \textbf{66} \\
        & SCIPsbt    & 1840 &  720.13 &  174.32 & 1793 &   68 &  380.43 &  179.95 &   \textbf{66} \\
        & SCIPabt    & 1840 &  \textbf{669.39} &  \textbf{169.55} & \textbf{1797} &   68 &  \textbf{373.90} &  \textbf{172.83} &   \textbf{66} \\
        \cmidrule{2-10}
        & GRBbasic   & 1840 &  601.86 &   \textbf{30.92} & 1816 &   68 &  \textbf{208.24} &   \textbf{15.38} &   68 \\
        & GRBsbt     & 1840 &  \textbf{315.93} &   70.03 & \textbf{1834} &   68 &  234.78 &   27.21 &   \textbf{66} \\
        \midrule
        \multirow{5}{*}{$s=4$}
        & SCIPbasic  & 1844 & 1745.56 &  487.75 & 1644 &   64 &  386.72 &  340.95 &   \textbf{64} \\
        & SCIPsbt    & 1844 & 1039.90 &  180.30 & 1734 &   64 &  155.26 &  144.16 &   \textbf{64} \\
        & SCIPabt    & 1844 & \textbf{1004.92} &  \textbf{177.84} & \textbf{1748} &   64 &  \textbf{147.46} &  \textbf{136.26} &   \textbf{64} \\
        \cmidrule{2-10}
        & GRBbasic   & 1844 & 1113.20 &   \textbf{58.17} & 1758 &   64 &   \textbf{11.90} &   \textbf{10.96} &   \textbf{64} \\
        & GRBsbt     & 1844 &  \textbf{420.78} &   84.82 & \textbf{1837} &   64 &   23.04 &   21.29 &   \textbf{64} \\
        \bottomrule
    \end{tabular*}
    \end{small}
    \end{center}
    \vskip -0.1in
\end{table*}

\begin{table}[t]
    \caption{Results for ENZYMES with local attack strength $s=2$ and different global budgets. The approaches tested are SCIP basic (SCIPbasic), SCIP static bounds tightening (SCIPsbt), SCIP aggressive bounds tightening (SCIPabt), Gurobi basic (GRBbasic), and Gurobi static bounds tightening (GRBsbt). For each global budget, the number of instances is 600, but we only present comparisons when all methods give consistent results except for time out, as shown in column ``\#". We compare times with respect to both average (``avg-time") and the shifted geometric mean (``sgm-time"), as well as the number of solved instances within time limits 2 hours (``\# solved"). Since robust instances are the ones where mixed-integer performance is most important (non-robust instances may frequently be found by more heuristic approaches), we also compare on the set of robust instances.}
    \vspace{-4mm}
    \label{tab:ENZYMES_s2_full}
    \vskip 0.15in
    \begin{center}
    \begin{scriptsize}
    \begin{tabular*}{\textwidth}{@{\;\;\extracolsep{\fill}}lrrrrrrrr}
    \toprule
    \multirow{2}{*}{method}
    & \multicolumn{4}{c}{all instances} & \multicolumn{4}{c}{robust instances}\\
    \cmidrule{2-5} \cmidrule{6-9}
     & \# & avg-time & sgm-time & \# solved & \# & avg-time & sgm-time & \# solved\\
    \midrule
    \multicolumn{9}{l}{global budget: $\lceil 0.01 \cdot |E| \rceil$}\\
    SCIPbasic  &  600 &  465.47 &   25.13 &  575 &  453 &  420.33 &   18.76 &  435 \\
    SCIPsbt    &  600 &  193.02 &   14.46 &  594 &  453 &  166.90 &   10.52 &  451 \\
    SCIPabt    &  600 &  190.05 &   14.58 &  592 &  453 &  161.85 &   10.53 &  449 \\
    GRBbasic   &  600 &  119.88 &    7.08 &  594 &  453 &   84.71 &    5.76 &  451 \\
    GRBsbt     &  600 &   87.31 &    5.90 &  596 &  453 &   47.40 &    4.89 &  453 \\
    \midrule
    \multicolumn{9}{l}{global budget: $\lceil 0.02 \cdot |E| \rceil$}\\
    SCIPbasic  &  597 &  630.47 &   37.42 &  557 &  392 &  457.99 &   18.92 &  371 \\
    SCIPsbt    &  597 &  248.67 &   20.65 &  585 &  392 &  181.41 &   10.44 &  385 \\
    SCIPabt    &  597 &  260.41 &   21.16 &  584 &  392 &  194.53 &   10.58 &  385 \\
    GRBbasic   &  597 &  170.91 &    9.36 &  588 &  392 &  111.73 &    5.05 &  388 \\
    GRBsbt     &  597 &  111.66 &    7.59 &  592 &  392 &   66.45 &    4.33 &  391 \\
    \midrule
    \multicolumn{9}{l}{global budget: $\lceil 0.03 \cdot |E| \rceil$}\\
    SCIPbasic  &  591 &  628.77 &   37.75 &  556 &  359 &  328.18 &   14.17 &  347 \\
    SCIPsbt    &  591 &  252.04 &   20.87 &  579 &  359 &  117.80 &    7.67 &  355 \\
    SCIPabt    &  591 &  262.35 &   21.03 &  577 &  359 &  120.84 &    7.67 &  354 \\
    GRBbasic   &  591 &  108.75 &    8.55 &  587 &  359 &   45.55 &    3.34 &  358 \\
    GRBsbt     &  591 &   79.77 &    7.05 &  588 &  359 &   22.37 &    2.79 &  359 \\
    \midrule
    \multicolumn{9}{l}{global budget: $\lceil 0.04 \cdot |E| \rceil$}\\
    SCIPbasic  &  593 &  633.12 &   39.29 &  557 &  346 &  256.36 &   11.96 &  336 \\
    SCIPsbt    &  593 &  234.00 &   21.84 &  585 &  346 &   61.20 &    6.07 &  344 \\
    SCIPabt    &  593 &  253.13 &   22.90 &  583 &  346 &   62.85 &    6.13 &  345 \\
    GRBbasic   &  593 &   80.87 &    7.59 &  591 &  346 &    6.44 &    2.19 &  346 \\
    GRBsbt     &  593 &   61.60 &    6.79 &  591 &  346 &    5.48 &    1.98 &  346 \\
    \midrule
    \multicolumn{9}{l}{global budget: $\lceil 0.05 \cdot |E| \rceil$}\\
    SCIPbasic  &  591 &  647.96 &   41.10 &  555 &  338 &  234.96 &   10.74 &  329 \\
    SCIPsbt    &  591 &  241.04 &   22.37 &  581 &  338 &   48.19 &    5.58 &  337 \\
    SCIPabt    &  591 &  293.33 &   23.46 &  576 &  338 &   75.63 &    5.66 &  335 \\
    GRBbasic   &  591 &   72.31 &    7.61 &  591 &  338 &    6.47 &    1.99 &  338 \\
    GRBsbt     &  591 &   70.81 &    7.34 &  589 &  338 &   13.86 &    1.83 &  338 \\
    \midrule
    \multicolumn{9}{l}{global budget: $\lceil 0.06 \cdot |E| \rceil$}\\
    SCIPbasic  &  594 &  623.67 &   41.01 &  560 &  334 &  204.73 &   10.09 &  327 \\
    SCIPsbt    &  594 &  210.26 &   22.99 &  588 &  334 &   37.46 &    5.22 &  333 \\
    SCIPabt    &  594 &  240.04 &   23.29 &  585 &  334 &   55.18 &    5.36 &  332 \\
    GRBbasic   &  594 &   73.46 &    7.56 &  594 &  334 &    8.41 &    1.85 &  334 \\
    GRBsbt     &  594 &   57.18 &    7.04 &  594 &  334 &   22.34 &    1.72 &  334 \\
    \midrule
    \multicolumn{9}{l}{global budget: $\lceil 0.07 \cdot |E| \rceil$}\\
    SCIPbasic  &  590 &  635.63 &   40.99 &  554 &  332 &  195.72 &    9.84 &  325 \\
    SCIPsbt    &  590 &  221.80 &   23.07 &  583 &  332 &   37.65 &    5.08 &  331 \\
    SCIPabt    &  590 &  259.83 &   23.68 &  583 &  332 &   54.13 &    5.25 &  330 \\
    GRBbasic   &  590 &   53.62 &    7.07 &  590 &  332 &    8.03 &    1.78 &  332 \\
    GRBsbt     &  590 &   62.36 &    7.46 &  590 &  332 &    6.69 &    1.77 &  332 \\
    \midrule
    \multicolumn{9}{l}{global budget: $\lceil 0.08 \cdot |E| \rceil$}\\
    SCIPbasic  &  583 &  608.75 &   39.97 &  553 &  332 &  197.92 &    9.95 &  325 \\
    SCIPsbt    &  583 &  240.35 &   22.52 &  575 &  332 &   38.82 &    5.11 &  331 \\
    SCIPabt    &  583 &  214.77 &   22.42 &  577 &  332 &   38.26 &    5.11 &  331 \\
    GRBbasic   &  583 &   60.73 &    6.92 &  582 &  332 &    8.32 &    1.85 &  332 \\
    GRBsbt     &  583 &   76.61 &    7.26 &  582 &  332 &    6.03 &    1.64 &  332 \\
    \midrule
    \multicolumn{9}{l}{global budget: $\lceil 0.09 \cdot |E| \rceil$}\\
    SCIPbasic  &  588 &  585.39 &   38.98 &  557 &  331 &  197.91 &   10.02 &  324 \\
    SCIPsbt    &  588 &  241.76 &   23.02 &  581 &  331 &   43.36 &    5.08 &  330 \\
    SCIPabt    &  588 &  249.56 &   23.67 &  580 &  331 &   37.99 &    5.05 &  330 \\
    GRBbasic   &  588 &   72.99 &    7.50 &  587 &  331 &    8.08 &    1.85 &  331 \\
    GRBsbt     &  588 &   75.35 &    7.23 &  586 &  331 &    6.08 &    1.66 &  331 \\
    \midrule
    \multicolumn{9}{l}{global budget: $\lceil 0.10 \cdot |E| \rceil$}\\
    SCIPbasic  &  588 &  602.12 &   39.14 &  555 &  332 &  205.53 &   10.01 &  325 \\
    SCIPsbt    &  588 &  223.54 &   22.05 &  580 &  332 &   46.58 &    5.09 &  331 \\
    SCIPabt    &  588 &  237.00 &   22.88 &  580 &  332 &   41.39 &    5.07 &  331 \\
    GRBbasic   &  588 &   44.81 &    6.83 &  588 &  332 &    7.72 &    1.96 &  332 \\
    GRBsbt     &  588 &   65.58 &    7.31 &  587 &  332 &   16.41 &    1.78 &  332 \\
    \bottomrule
  \end{tabular*}
  \end{scriptsize}
  \end{center}
  \vskip -0.1in
\end{table}

\begin{table}[t]
    \caption{Results for ENZYMES with local attack strength $s=3$ and different global budgets. The approaches tested are SCIP basic (SCIPbasic), SCIP static bounds tightening (SCIPsbt), SCIP aggressive bounds tightening (SCIPabt), Gurobi basic (GRBbasic), and Gurobi static bounds tightening (GRBsbt). For each global budget, the number of instances is 600, but we only present comparisons when all methods give consistent results except for time out, as shown in column ``\#". We compare times with respect to both average (``avg-time") and the shifted geometric mean (``sgm-time"), as well as the number of solved instances within time limits 2 hours (``\# solved"). Since robust instances are the ones where mixed-integer performance is most important (non-robust instances may frequently be found by more heuristic approaches), we also compare on the set of robust instances.}
    \vspace{-4mm}
    \label{tab:ENZYMES_s3_full}
    \vskip 0.15in
    \begin{center}
    \begin{scriptsize}
    \begin{tabular*}{\textwidth}{@{\;\;\extracolsep{\fill}}lrrrrrrrr}
    \toprule
    \multirow{2}{*}{method}
    & \multicolumn{4}{c}{all instances} & \multicolumn{4}{c}{robust instances}\\
    \cmidrule{2-5} \cmidrule{6-9}
     & \# & avg-time & sgm-time & \# solved & \# & avg-time & sgm-time & \# solved\\
    \midrule
    \multicolumn{9}{l}{global budget: $\lceil 0.01 \cdot |E| \rceil$}\\
    SCIPbasic  &  597 & 2041.24 &  248.53 &  464 &  319 & 1088.49 &   94.41 &  287 \\
    SCIPsbt    &  597 & 1160.38 &  109.23 &  521 &  319 &  481.74 &   43.92 &  307 \\
    SCIPabt    &  597 & 1129.15 &  109.77 &  523 &  319 &  478.08 &   43.68 &  307 \\
    GRBbasic   &  597 &  985.15 &   52.53 &  542 &  319 &  370.09 &   20.68 &  310 \\
    GRBsbt     &  597 &  700.88 &   39.96 &  558 &  319 &  193.68 &   16.28 &  318 \\
    \midrule
    \multicolumn{9}{l}{global budget: $\lceil 0.02 \cdot |E| \rceil$}\\
    SCIPbasic  &  587 & 2738.05 &  498.94 &  413 &  218 & 1896.50 &  185.55 &  179 \\
    SCIPsbt    &  587 & 1683.98 &  226.92 &  484 &  218 & 1004.56 &   88.59 &  199 \\
    SCIPabt    &  587 & 1623.45 &  227.82 &  491 &  218 &  970.46 &   86.96 &  207 \\
    GRBbasic   &  587 & 1059.95 &   99.76 &  544 &  218 &  495.46 &   43.64 &  215 \\
    GRBsbt     &  587 &  728.29 &   66.57 &  552 &  218 &  352.89 &   34.88 &  217 \\
    \midrule
    \multicolumn{9}{l}{global budget: $\lceil 0.03 \cdot |E| \rceil$}\\
    SCIPbasic  &  590 & 2835.26 &  542.92 &  403 &  169 & 1685.32 &  142.10 &  139 \\
    SCIPsbt    &  590 & 1745.83 &  245.95 &  482 &  169 &  883.20 &   67.97 &  157 \\
    SCIPabt    &  590 & 1757.43 &  252.98 &  477 &  169 &  857.83 &   66.17 &  155 \\
    GRBbasic   &  590 &  926.16 &   94.63 &  551 &  169 &  346.30 &   30.00 &  166 \\
    GRBsbt     &  590 &  597.79 &   64.60 &  565 &  169 &  266.98 &   25.74 &  168 \\
    \midrule
    \multicolumn{9}{l}{global budget: $\lceil 0.04 \cdot |E| \rceil$}\\
    SCIPbasic  &  588 & 2793.58 &  539.45 &  404 &  141 & 1623.34 &  113.29 &  114 \\
    SCIPsbt    &  588 & 1684.95 &  241.19 &  486 &  141 &  844.50 &   54.50 &  129 \\
    SCIPabt    &  588 & 1779.51 &  259.53 &  481 &  141 &  869.72 &   53.91 &  128 \\
    GRBbasic   &  588 &  705.60 &   84.77 &  568 &  141 &  302.69 &   23.50 &  141 \\
    GRBsbt     &  588 &  498.54 &   63.82 &  572 &  141 &  268.36 &   20.55 &  139 \\
    \midrule
    \multicolumn{9}{l}{global budget: $\lceil 0.05 \cdot |E| \rceil$}\\
    SCIPbasic  &  589 & 2730.97 &  502.64 &  407 &  128 & 1462.98 &   88.07 &  104 \\
    SCIPsbt    &  589 & 1559.45 &  230.58 &  499 &  128 &  825.88 &   44.16 &  117 \\
    SCIPabt    &  589 & 1588.79 &  241.82 &  500 &  128 &  821.32 &   43.98 &  116 \\
    GRBbasic   &  589 &  653.18 &   82.60 &  572 &  128 &  363.43 &   20.28 &  128 \\
    GRBsbt     &  589 &  597.25 &   69.81 &  571 &  128 &  282.16 &   17.86 &  126 \\
    \midrule
    \multicolumn{9}{l}{global budget: $\lceil 0.06 \cdot |E| \rceil$}\\
    SCIPbasic  &  587 & 2678.33 &  483.01 &  409 &  120 & 1405.21 &   82.54 &   99 \\
    SCIPsbt    &  587 & 1523.67 &  226.15 &  502 &  120 &  798.68 &   41.36 &  110 \\
    SCIPabt    &  587 & 1611.51 &  243.19 &  496 &  120 &  813.46 &   42.08 &  110 \\
    GRBbasic   &  587 &  719.60 &   81.19 &  564 &  120 &  415.24 &   19.88 &  118 \\
    GRBsbt     &  587 &  664.13 &   78.00 &  558 &  120 &  214.15 &   17.07 &  120 \\
    \midrule
    \multicolumn{9}{l}{global budget: $\lceil 0.07 \cdot |E| \rceil$}\\
    SCIPbasic  &  576 & 2729.50 &  492.39 &  400 &  115 & 1290.92 &   70.80 &   97 \\
    SCIPsbt    &  576 & 1518.06 &  226.43 &  493 &  115 &  699.68 &   35.11 &  105 \\
    SCIPabt    &  576 & 1605.08 &  240.04 &  489 &  115 &  682.72 &   34.36 &  106 \\
    GRBbasic   &  576 &  762.36 &   88.48 &  552 &  115 &  332.71 &   16.17 &  114 \\
    GRBsbt     &  576 &  780.03 &   86.20 &  541 &  115 &  209.06 &   13.87 &  114 \\
    \midrule
    \multicolumn{9}{l}{global budget: $\lceil 0.08 \cdot |E| \rceil$}\\
    SCIPbasic  &  579 & 2527.15 &  450.55 &  425 &  116 & 1288.86 &   70.91 &   97 \\
    SCIPsbt    &  579 & 1445.35 &  212.99 &  499 &  116 &  675.81 &   35.33 &  107 \\
    SCIPabt    &  579 & 1492.87 &  226.64 &  500 &  116 &  700.39 &   36.03 &  107 \\
    GRBbasic   &  579 &  828.08 &   90.22 &  552 &  116 &  307.12 &   16.05 &  114 \\
    GRBsbt     &  579 &  842.19 &   93.20 &  539 &  116 &  146.60 &   13.40 &  116 \\
    \midrule
    \multicolumn{9}{l}{global budget: $\lceil 0.09 \cdot |E| \rceil$}\\
    SCIPbasic  &  586 & 2594.38 &  476.06 &  415 &  113 & 1190.02 &   64.52 &   96 \\
    SCIPsbt    &  586 & 1498.64 &  232.34 &  512 &  113 &  701.85 &   33.91 &  104 \\
    SCIPabt    &  586 & 1634.99 &  246.41 &  499 &  113 &  620.21 &   32.32 &  106 \\
    GRBbasic   &  586 &  863.88 &   96.14 &  561 &  113 &  260.83 &   14.23 &  112 \\
    GRBsbt     &  586 &  925.87 &   99.95 &  539 &  113 &  178.31 &   12.64 &  113 \\
    \midrule
    \multicolumn{9}{l}{global budget: $\lceil 0.10 \cdot |E| \rceil$}\\
    SCIPbasic  &  576 & 2627.36 &  459.95 &  409 &  115 & 1292.84 &   69.60 &   96 \\
    SCIPsbt    &  576 & 1513.94 &  234.73 &  499 &  115 &  746.69 &   35.91 &  105 \\
    SCIPabt    &  576 & 1620.34 &  252.07 &  487 &  115 &  727.18 &   35.21 &  105 \\
    GRBbasic   &  576 &  952.10 &   99.60 &  543 &  115 &  388.85 &   16.03 &  112 \\
    GRBsbt     &  576 & 1058.39 &  112.22 &  529 &  115 &  225.30 &   13.55 &  114 \\
    \bottomrule
  \end{tabular*}
  \end{scriptsize}
  \end{center}
  \vskip -0.1in
\end{table}

\begin{table}[t]
    \caption{Results for ENZYMES with local attack strength $s=4$ and different global budgets. The approaches tested are SCIP basic (SCIPbasic), SCIP static bounds tightening (SCIPsbt), SCIP aggressive bounds tightening (SCIPabt), Gurobi basic (GRBbasic), and Gurobi static bounds tightening (GRBsbt). For each global budget, the number of instances is 600, but we only present comparisons when all methods give consistent results except for time out, as shown in column ``\#". We compare times with respect to both average (``avg-time") and the shifted geometric mean (``sgm-time"), as well as the number of solved instances within time limits 2 hours (``\# solved"). Since robust instances are the ones where mixed-integer performance is most important (non-robust instances may frequently be found by more heuristic approaches), we also compare on the set of robust instances.}
    \vspace{-4mm}
    \label{tab:ENZYMES_s4_full}
    \vskip 0.15in
    \begin{center}
    \begin{scriptsize}
    \begin{tabular*}{\textwidth}{@{\;\;\extracolsep{\fill}}lrrrrrrrr}
    \toprule
    \multirow{2}{*}{method}
    & \multicolumn{4}{c}{all instances} & \multicolumn{4}{c}{robust instances}\\
    \cmidrule{2-5} \cmidrule{6-9}
     & \# & avg-time & sgm-time & \# solved & \# & avg-time & sgm-time & \# solved\\
    \midrule
    \multicolumn{9}{l}{global budget: $\lceil 0.01 \cdot |E| \rceil$}\\
    SCIPbasic  &  594 & 3421.76 &  773.05 &  345 &  211 & 1260.49 &  170.39 &  182 \\
    SCIPsbt    &  594 & 2352.82 &  345.23 &  425 &  211 &  593.60 &   74.72 &  202 \\
    SCIPabt    &  594 & 2282.73 &  337.01 &  432 &  211 &  569.32 &   74.28 &  203 \\
    GRBbasic   &  594 & 2107.45 &  170.88 &  456 &  211 &  470.19 &   24.40 &  204 \\
    GRBsbt     &  594 & 1595.10 &  130.26 &  493 &  211 &  260.64 &   20.56 &  211 \\
    \midrule
    \multicolumn{9}{l}{global budget: $\lceil 0.02 \cdot |E| \rceil$}\\
    SCIPbasic  &  593 & 4550.98 & 1874.90 &  280 &  103 & 2833.99 &  581.49 &   76 \\
    SCIPsbt    &  593 & 3259.56 &  770.71 &  362 &  103 & 2135.50 &  324.52 &   82 \\
    SCIPabt    &  593 & 3120.26 &  741.51 &  374 &  103 & 2049.08 &  316.62 &   84 \\
    GRBbasic   &  593 & 2373.24 &  422.28 &  463 &  103 & 1195.25 &  153.67 &   98 \\
    GRBsbt     &  593 & 1658.01 &  233.42 &  502 &  103 &  819.20 &  108.68 &  102 \\
    \midrule
    \multicolumn{9}{l}{global budget: $\lceil 0.03 \cdot |E| \rceil$}\\
    SCIPbasic  &  591 & 4476.55 & 1883.37 &  270 &   65 & 2733.32 &  583.01 &   45 \\
    SCIPsbt    &  591 & 3307.01 &  869.32 &  360 &   65 & 1900.42 &  315.89 &   54 \\
    SCIPabt    &  591 & 3340.71 &  888.09 &  353 &   65 & 1820.57 &  298.76 &   54 \\
    GRBbasic   &  591 & 1984.46 &  403.78 &  504 &   65 & 1329.91 &  164.43 &   60 \\
    GRBsbt     &  591 & 1283.41 &  191.51 &  531 &   65 &  993.78 &  123.35 &   63 \\
    \midrule
    \multicolumn{9}{l}{global budget: $\lceil 0.04 \cdot |E| \rceil$}\\
    SCIPbasic  &  593 & 4469.33 & 1890.99 &  268 &   46 & 3027.57 &  602.22 &   28 \\
    SCIPsbt    &  593 & 3339.44 &  808.49 &  357 &   46 & 2376.84 &  332.05 &   35 \\
    SCIPabt    &  593 & 3312.35 &  824.98 &  362 &   46 & 2295.64 &  327.83 &   36 \\
    GRBbasic   &  593 & 1643.02 &  335.49 &  537 &   46 & 1407.86 &  154.94 &   41 \\
    GRBsbt     &  593 & 1039.69 &  175.10 &  550 &   46 & 1024.66 &  106.28 &   44 \\
    \midrule
    \multicolumn{9}{l}{global budget: $\lceil 0.05 \cdot |E| \rceil$}\\
    SCIPbasic  &  595 & 4379.75 & 1691.84 &  286 &   34 & 2571.61 &  438.08 &   23 \\
    SCIPsbt    &  595 & 2956.30 &  659.83 &  396 &   34 & 1625.51 &  205.03 &   30 \\
    SCIPabt    &  595 & 3034.27 &  691.12 &  385 &   34 & 1651.50 &  201.25 &   29 \\
    GRBbasic   &  595 & 1538.31 &  300.14 &  539 &   34 &  843.47 &   86.46 &   32 \\
    GRBsbt     &  595 & 1035.63 &  198.01 &  561 &   34 &  696.91 &   67.81 &   34 \\
    \midrule
    \multicolumn{9}{l}{global budget: $\lceil 0.06 \cdot |E| \rceil$}\\
    SCIPbasic  &  588 & 4108.15 & 1427.47 &  302 &   29 & 2220.69 &  366.12 &   21 \\
    SCIPsbt    &  588 & 2879.51 &  624.42 &  399 &   29 & 1571.68 &  174.19 &   24 \\
    SCIPabt    &  588 & 2891.41 &  644.29 &  404 &   29 & 1619.69 &  174.26 &   25 \\
    GRBbasic   &  588 & 1472.87 &  283.85 &  543 &   29 & 1209.09 &   79.81 &   27 \\
    GRBsbt     &  588 & 1035.34 &  214.00 &  551 &   29 &  518.49 &   49.90 &   29 \\
    \midrule
    \multicolumn{9}{l}{global budget: $\lceil 0.07 \cdot |E| \rceil$}\\
    SCIPbasic  &  592 & 4147.98 & 1484.03 &  309 &   24 & 1749.90 &  288.38 &   19 \\
    SCIPsbt    &  592 & 2685.22 &  543.80 &  414 &   24 & 1005.48 &  122.89 &   22 \\
    SCIPabt    &  592 & 2732.23 &  556.59 &  416 &   24 & 1053.82 &  119.60 &   22 \\
    GRBbasic   &  592 & 1609.26 &  294.58 &  532 &   24 &  679.49 &   41.95 &   24 \\
    GRBsbt     &  592 & 1298.98 &  257.15 &  541 &   24 &  155.27 &   32.76 &   24 \\
    \midrule
    \multicolumn{9}{l}{global budget: $\lceil 0.08 \cdot |E| \rceil$}\\
    SCIPbasic  &  588 & 4085.50 & 1440.95 &  308 &   22 & 1622.94 &  281.69 &   18 \\
    SCIPsbt    &  588 & 2489.10 &  505.68 &  430 &   22 & 1091.16 &  114.08 &   19 \\
    SCIPabt    &  588 & 2693.02 &  563.16 &  418 &   22 & 1089.48 &  114.84 &   19 \\
    GRBbasic   &  588 & 1574.46 &  297.33 &  532 &   22 &  472.77 &   33.43 &   21 \\
    GRBsbt     &  588 & 1399.73 &  308.64 &  536 &   22 &  502.12 &   34.93 &   22 \\
    \midrule
    \multicolumn{9}{l}{global budget: $\lceil 0.09 \cdot |E| \rceil$}\\
    SCIPbasic  &  579 & 4023.93 & 1476.52 &  310 &   22 & 1739.77 &  287.68 &   18 \\
    SCIPsbt    &  579 & 2443.85 &  480.58 &  422 &   22 & 1384.42 &  127.34 &   18 \\
    SCIPabt    &  579 & 2599.45 &  530.96 &  427 &   22 & 1384.49 &  127.40 &   18 \\
    GRBbasic   &  579 & 1726.07 &  315.56 &  504 &   22 &  920.76 &   45.09 &   20 \\
    GRBsbt     &  579 & 1574.30 &  346.97 &  517 &   22 &  961.86 &   45.76 &   21 \\
    \midrule
    \multicolumn{9}{l}{global budget: $\lceil 0.10 \cdot |E| \rceil$}\\
    SCIPbasic  &  588 & 3971.18 & 1391.04 &  325 &   21 & 1375.79 &  221.98 &   18 \\
    SCIPsbt    &  588 & 2294.64 &  426.34 &  447 &   21 & 1111.89 &  107.09 &   18 \\
    SCIPabt    &  588 & 2398.19 &  460.41 &  442 &   21 & 1079.78 &  103.04 &   19 \\
    GRBbasic   &  588 & 1772.66 &  311.55 &  515 &   21 &  632.60 &   34.92 &   20 \\
    GRBsbt     &  588 & 1814.80 &  403.30 &  508 &   21 &  695.35 &   34.77 &   20 \\
    \bottomrule
  \end{tabular*}
  \end{scriptsize}
  \end{center}
  \vskip -0.1in
\end{table}

\begin{table}[t]
    \caption{Results for MUTAG with local attack strength $s=2$ and different global budgets. The approaches tested are SCIP basic (SCIPbasic), SCIP static bounds tightening (SCIPsbt), SCIP aggressive bounds tightening (SCIPabt), Gurobi basic (GRBbasic), and Gurobi static bounds tightening (GRBsbt). For each global budget, the number of instances is 188, but we only present comparisons when all methods give consistent results except for time out, as shown in column ``\#". We compare times with respect to both average (``avg-time") and the shifted geometric mean (``sgm-time"), as well as the number of solved instances within time limits 2 hours (``\# solved"). Since robust instances are the ones where mixed-integer performance is most important (non-robust instances may frequently be found by more heuristic approaches), we also compare on the set of robust instances.}
    \vspace{-4mm}
    \label{tab:MUTAG_s2_full}
    \vskip 0.15in
    \begin{center}
    \begin{scriptsize}
    \begin{tabular*}{\textwidth}{@{\;\;\extracolsep{\fill}}lrrrrrrrr}
    \toprule
    \multirow{2}{*}{method}
    & \multicolumn{4}{c}{all instances} & \multicolumn{4}{c}{robust instances}\\
    \cmidrule{2-5} \cmidrule{6-9}
    & \# & avg-time & sgm-time & \# solved & \# & avg-time & sgm-time & \# solved\\
    \midrule
    \multicolumn{9}{l}{global budget: $\lceil 0.01 \cdot |E| \rceil$}\\
    SCIPbasic  &  124 &   74.28 &   45.27 &  124 &   24 &  176.24 &  162.61 &   24 \\
    SCIPsbt    &  124 &   33.89 &   22.76 &  124 &   24 &   84.37 &   79.12 &   24 \\
    SCIPabt    &  124 &   32.31 &   22.18 &  124 &   24 &   77.22 &   72.94 &   24 \\
    GRBbasic   &  124 &    2.30 &    2.16 &  124 &   24 &    4.91 &    4.82 &   24 \\
    GRBsbt     &  124 &    4.10 &    3.69 &  124 &   24 &    9.39 &    8.98 &   24 \\
    \midrule
    \multicolumn{9}{l}{global budget: $\lceil 0.02 \cdot |E| \rceil$}\\
    SCIPbasic  &  141 &  437.96 &   87.20 &  139 &   16 & 1028.39 &  226.02 &   14 \\
    SCIPsbt    &  141 &  127.08 &   37.84 &  141 &   16 &  407.73 &  115.97 &   16 \\
    SCIPabt    &  141 &   98.31 &   36.72 &  141 &   16 &  287.84 &  103.45 &   16 \\
    GRBbasic   &  141 &   94.67 &   11.83 &  141 &   16 &  205.18 &   15.21 &   16 \\
    GRBsbt     &  141 &   26.14 &    8.90 &  141 &   16 &   93.51 &   17.19 &   16 \\
    \midrule
    \multicolumn{9}{l}{global budget: $\lceil 0.03 \cdot |E| \rceil$}\\
    SCIPbasic  &  180 &  542.88 &  162.56 &  178 &    3 & 4807.52 & 1181.75 &    1 \\
    SCIPsbt    &  180 &  152.61 &   57.72 &  180 &    3 & 2080.69 &  545.97 &    3 \\
    SCIPabt    &  180 &  134.96 &   58.51 &  180 &    3 & 1334.79 &  425.77 &    3 \\
    GRBbasic   &  180 &   87.25 &   10.17 &  180 &    3 & 1245.54 &  325.46 &    3 \\
    GRBsbt     &  180 &   34.48 &   16.03 &  180 &    3 &  507.82 &  176.89 &    3 \\
    \midrule
    \multicolumn{9}{l}{global budget: $\lceil 0.04 \cdot |E| \rceil$}\\
    SCIPbasic  &  174 &  613.49 &  199.84 &  172 &    1 &   26.09 &   26.09 &    1 \\
    SCIPsbt    &  174 &  175.59 &   66.82 &  174 &    1 &   11.45 &   11.45 &    1 \\
    SCIPabt    &  174 &  133.90 &   67.08 &  174 &    1 &   13.23 &   13.23 &    1 \\
    GRBbasic   &  174 &   49.70 &    4.59 &  174 &    1 &    1.51 &    1.51 &    1 \\
    GRBsbt     &  174 &   47.76 &   21.59 &  174 &    1 &    1.55 &    1.55 &    1 \\
    \midrule
    \multicolumn{9}{l}{global budget: $\lceil 0.05 \cdot |E| \rceil$}\\
    SCIPbasic  &  167 &  707.14 &  233.17 &  166 &    0 & ---     & ---     & ---  \\
    SCIPsbt    &  167 &  185.08 &   81.61 &  167 &    0 & ---     & ---     & ---  \\
    SCIPabt    &  167 &  170.37 &   82.49 &  167 &    0 & ---     & ---     & ---  \\
    GRBbasic   &  167 &   23.07 &    2.30 &  167 &    0 & ---     & ---     & ---  \\
    GRBsbt     &  167 &   61.18 &   27.22 &  167 &    0 & ---     & ---     & ---  \\
    \midrule
    \multicolumn{9}{l}{global budget: $\lceil 0.06 \cdot |E| \rceil$}\\
    SCIPbasic  &  158 &  832.55 &  270.16 &  157 &    0 & ---     & ---     & ---  \\
    SCIPsbt    &  158 &  190.53 &   92.80 &  158 &    0 & ---     & ---     & ---  \\
    SCIPabt    &  158 &  215.46 &  102.23 &  158 &    0 & ---     & ---     & ---  \\
    GRBbasic   &  158 &   23.74 &    2.55 &  158 &    0 & ---     & ---     & ---  \\
    GRBsbt     &  158 &   57.60 &   27.95 &  158 &    0 & ---     & ---     & ---  \\
    \midrule
    \multicolumn{9}{l}{global budget: $\lceil 0.07 \cdot |E| \rceil$}\\
    SCIPbasic  &  159 &  841.40 &  265.97 &  158 &    0 & ---     & ---     & ---  \\
    SCIPsbt    &  159 &  219.35 &   94.60 &  159 &    0 & ---     & ---     & ---  \\
    SCIPabt    &  159 &  264.65 &  112.42 &  159 &    0 & ---     & ---     & ---  \\
    GRBbasic   &  159 &   43.72 &    2.98 &  159 &    0 & ---     & ---     & ---  \\
    GRBsbt     &  159 &   57.12 &   25.68 &  159 &    0 & ---     & ---     & ---  \\
    \midrule
    \multicolumn{9}{l}{global budget: $\lceil 0.08 \cdot |E| \rceil$}\\
    SCIPbasic  &  159 &  859.84 &  265.24 &  159 &    0 & ---     & ---     & ---  \\
    SCIPsbt    &  159 &  254.42 &  109.18 &  159 &    0 & ---     & ---     & ---  \\
    SCIPabt    &  159 &  299.46 &  134.91 &  159 &    0 & ---     & ---     & ---  \\
    GRBbasic   &  159 &    4.31 &    2.18 &  159 &    0 & ---     & ---     & ---  \\
    GRBsbt     &  159 &   65.85 &   27.45 &  159 &    0 & ---     & ---     & ---  \\
    \midrule
    \multicolumn{9}{l}{global budget: $\lceil 0.09 \cdot |E| \rceil$}\\
    SCIPbasic  &  165 &  877.25 &  254.34 &  163 &    0 & ---     & ---     & ---  \\
    SCIPsbt    &  165 &  274.26 &  113.80 &  165 &    0 & ---     & ---     & ---  \\
    SCIPabt    &  165 &  333.40 &  137.86 &  165 &    0 & ---     & ---     & ---  \\
    GRBbasic   &  165 &    4.26 &    2.07 &  165 &    0 & ---     & ---     & ---  \\
    GRBsbt     &  165 &   89.31 &   31.92 &  165 &    0 & ---     & ---     & ---  \\
    \midrule
    \multicolumn{9}{l}{global budget: $\lceil 0.10 \cdot |E| \rceil$}\\
    SCIPbasic  &  162 &  864.13 &  251.70 &  159 &    0 & ---     & ---     & ---  \\
    SCIPsbt    &  162 &  307.52 &  120.64 &  162 &    0 & ---     & ---     & ---  \\
    SCIPabt    &  162 &  352.19 &  149.21 &  162 &    0 & ---     & ---     & ---  \\
    GRBbasic   &  162 &    6.29 &    2.89 &  162 &    0 & ---     & ---     & ---  \\
    GRBsbt     &  162 &  141.44 &   49.11 &  162 &    0 & ---     & ---     & ---  \\
    \bottomrule
  \end{tabular*}
  \end{scriptsize}
  \end{center}
  \vskip -0.1in
\end{table}  

\begin{table}[t]
    \caption{Results for MUTAG with local attack strength $s=3$ and different global budgets. The approaches tested are SCIP basic (SCIPbasic), SCIP static bounds tightening (SCIPsbt), SCIP aggressive bounds tightening (SCIPabt), Gurobi basic (GRBbasic), and Gurobi static bounds tightening (GRBsbt). For each global budget, the number of instances is 188, but we only present comparisons when all methods give consistent results except for time out, as shown in column ``\#". We compare times with respect to both average (``avg-time") and the shifted geometric mean (``sgm-time"), as well as the number of solved instances within time limits 2 hours (``\# solved"). Since robust instances are the ones where mixed-integer performance is most important (non-robust instances may frequently be found by more heuristic approaches), we also compare on the set of robust instances.}
    \vspace{-4mm}
    \label{tab:MUTAG_s3_full}
    \vskip 0.15in
    \begin{center}
    \begin{scriptsize}
    \begin{tabular*}{\textwidth}{@{\;\;\extracolsep{\fill}}lrrrrrrrr}
    \toprule
    \multirow{2}{*}{method}
    & \multicolumn{4}{c}{all instances} & \multicolumn{4}{c}{robust instances}\\
    \cmidrule{2-5} \cmidrule{6-9}
    & \# & avg-time & sgm-time & \# solved & \# & avg-time & sgm-time & \# solved\\
    \midrule
    \multicolumn{9}{l}{global budget: $\lceil 0.01 \cdot |E| \rceil$}\\
    SCIPbasic  &  178 &  216.94 &  138.74 &  178 &   43 &  484.49 &  438.14 &   43 \\
    SCIPsbt    &  178 &   97.37 &   68.28 &  178 &   43 &  191.70 &  180.49 &   43 \\
    SCIPabt    &  178 &   96.24 &   67.91 &  178 &   43 &  186.38 &  174.65 &   43 \\
    GRBbasic   &  178 &    6.75 &    5.92 &  178 &   43 &   13.62 &   12.68 &   43 \\
    GRBsbt     &  178 &   11.81 &    9.64 &  178 &   43 &   25.98 &   24.19 &   43 \\
    \midrule
    \multicolumn{9}{l}{global budget: $\lceil 0.02 \cdot |E| \rceil$}\\
    SCIPbasic  &  183 &  964.93 &  197.33 &  169 &   22 &  709.48 &  418.88 &   21 \\
    SCIPsbt    &  183 &  353.80 &   93.62 &  180 &   22 &  470.79 &  172.81 &   21 \\
    SCIPabt    &  183 &  348.24 &   92.84 &  180 &   22 &  461.02 &  162.62 &   21 \\
    GRBbasic   &  183 &  494.45 &   23.69 &  181 &   22 &  315.86 &   16.25 &   22 \\
    GRBsbt     &  183 &  216.83 &   20.71 &  180 &   22 &  347.31 &   27.91 &   21 \\
    \midrule
    \multicolumn{9}{l}{global budget: $\lceil 0.03 \cdot |E| \rceil$}\\
    SCIPbasic  &  185 & 1457.00 &  260.71 &  169 &    2 & 3634.79 &  747.48 &    1 \\
    SCIPsbt    &  185 &  424.80 &  119.43 &  182 &    2 & 3616.68 &  549.13 &    1 \\
    SCIPabt    &  185 &  418.05 &  118.39 &  182 &    2 & 3616.66 &  548.81 &    1 \\
    GRBbasic   &  185 &  744.52 &   83.02 &  183 &    2 & 3311.59 &  274.89 &    2 \\
    GRBsbt     &  185 &  257.34 &   44.27 &  182 &    2 & 3601.72 &  301.41 &    1 \\
    \midrule
    \multicolumn{9}{l}{global budget: $\lceil 0.04 \cdot |E| \rceil$}\\
    SCIPbasic  &  185 & 1473.13 &  290.37 &  168 &    1 &   73.10 &   73.10 &    1 \\
    SCIPsbt    &  185 &  605.01 &  151.85 &  179 &    1 &   35.09 &   35.09 &    1 \\
    SCIPabt    &  185 &  557.49 &  149.98 &  182 &    1 &   35.24 &   35.24 &    1 \\
    GRBbasic   &  185 &  910.40 &   82.52 &  182 &    1 &    2.25 &    2.25 &    1 \\
    GRBsbt     &  185 &  204.63 &   52.12 &  185 &    1 &    3.55 &    3.55 &    1 \\
    \midrule
    \multicolumn{9}{l}{global budget: $\lceil 0.05 \cdot |E| \rceil$}\\
    SCIPbasic  &  186 & 1464.64 &  307.77 &  167 &    0 & ---     & ---     & ---  \\
    SCIPsbt    &  186 &  661.42 &  211.00 &  182 &    0 & ---     & ---     & ---  \\
    SCIPabt    &  186 &  637.86 &  204.26 &  183 &    0 & ---     & ---     & ---  \\
    GRBbasic   &  186 & 1072.39 &   91.48 &  183 &    0 & ---     & ---     & ---  \\
    GRBsbt     &  186 &  282.76 &   77.28 &  186 &    0 & ---     & ---     & ---  \\
    \midrule
    \multicolumn{9}{l}{global budget: $\lceil 0.06 \cdot |E| \rceil$}\\
    SCIPbasic  &  186 & 1341.58 &  261.62 &  171 &    0 & ---     & ---     & ---  \\
    SCIPsbt    &  186 &  700.35 &  212.61 &  181 &    0 & ---     & ---     & ---  \\
    SCIPabt    &  186 &  626.75 &  198.82 &  183 &    0 & ---     & ---     & ---  \\
    GRBbasic   &  186 &  934.94 &   59.53 &  184 &    0 & ---     & ---     & ---  \\
    GRBsbt     &  186 &  334.13 &   93.00 &  186 &    0 & ---     & ---     & ---  \\
    \midrule
    \multicolumn{9}{l}{global budget: $\lceil 0.07 \cdot |E| \rceil$}\\
    SCIPbasic  &  186 & 1537.24 &  306.34 &  166 &    0 & ---     & ---     & ---  \\
    SCIPsbt    &  186 &  964.08 &  244.13 &  181 &    0 & ---     & ---     & ---  \\
    SCIPabt    &  186 &  919.51 &  235.79 &  180 &    0 & ---     & ---     & ---  \\
    GRBbasic   &  186 &  564.47 &   23.72 &  183 &    0 & ---     & ---     & ---  \\
    GRBsbt     &  186 &  375.94 &  104.92 &  186 &    0 & ---     & ---     & ---  \\
    \midrule
    \multicolumn{9}{l}{global budget: $\lceil 0.08 \cdot |E| \rceil$}\\
    SCIPbasic  &  186 & 1559.03 &  279.39 &  166 &    0 & ---     & ---     & ---  \\
    SCIPsbt    &  186 & 1065.84 &  242.82 &  179 &    0 & ---     & ---     & ---  \\
    SCIPabt    &  186 &  989.15 &  236.15 &  179 &    0 & ---     & ---     & ---  \\
    GRBbasic   &  186 &  388.40 &   13.63 &  184 &    0 & ---     & ---     & ---  \\
    GRBsbt     &  186 &  426.44 &  132.03 &  186 &    0 & ---     & ---     & ---  \\
    \midrule
    \multicolumn{9}{l}{global budget: $\lceil 0.09 \cdot |E| \rceil$}\\
    SCIPbasic  &  184 & 1720.80 &  331.98 &  161 &    0 & ---     & ---     & ---  \\
    SCIPsbt    &  184 & 1080.19 &  279.90 &  179 &    0 & ---     & ---     & ---  \\
    SCIPabt    &  184 &  993.48 &  271.16 &  178 &    0 & ---     & ---     & ---  \\
    GRBbasic   &  184 &  456.38 &   12.95 &  182 &    0 & ---     & ---     & ---  \\
    GRBsbt     &  184 &  496.92 &  164.96 &  184 &    0 & ---     & ---     & ---  \\
    \midrule
    \multicolumn{9}{l}{global budget: $\lceil 0.10 \cdot |E| \rceil$}\\
    SCIPbasic  &  181 & 1777.61 &  317.24 &  155 &    0 & ---     & ---     & ---  \\
    SCIPsbt    &  181 & 1231.10 &  271.38 &  172 &    0 & ---     & ---     & ---  \\
    SCIPabt    &  181 & 1090.15 &  259.43 &  172 &    0 & ---     & ---     & ---  \\
    GRBbasic   &  181 &  414.44 &   10.81 &  176 &    0 & ---     & ---     & ---  \\
    GRBsbt     &  181 &  545.06 &  187.13 &  181 &    0 & ---     & ---     & ---  \\
    \bottomrule
  \end{tabular*}
  \end{scriptsize}
  \end{center}
  \vskip -0.1in
\end{table}  

\begin{table}[t]
    \caption{Results for MUTAG with local attack strength $s=4$ and different global budgets. The approaches tested are SCIP basic (SCIPbasic), SCIP static bounds tightening (SCIPsbt), SCIP aggressive bounds tightening (SCIPabt), Gurobi basic (GRBbasic), and Gurobi static bounds tightening (GRBsbt). For each global budget, the number of instances is 188, but we only present comparisons when all methods give consistent results except for time out, as shown in column ``\#". We compare times with respect to both average (``avg-time") and the shifted geometric mean (``sgm-time"), as well as the number of solved instances within time limits 2 hours (``\# solved"). Since robust instances are the ones where mixed-integer performance is most important (non-robust instances may frequently be found by more heuristic approaches), we also compare on the set of robust instances.}
    \vspace{-4mm}
    \label{tab:MUTAG_s4_full}
    \vskip 0.15in
    \begin{center}
    \begin{scriptsize}
    \begin{tabular*}{\textwidth}{@{\;\;\extracolsep{\fill}}lrrrrrrrr}
    \toprule
    \multirow{2}{*}{method}
    & \multicolumn{4}{c}{all instances} & \multicolumn{4}{c}{robust instances}\\
    \cmidrule{2-5} \cmidrule{6-9}
    & \# & avg-time & sgm-time & \# solved & \# & avg-time & sgm-time & \# solved\\
    \midrule
    \multicolumn{9}{l}{global budget: $\lceil 0.01 \cdot |E| \rceil$}\\
    SCIPbasic  &  178 &  191.35 &  123.98 &  178 &   42 &  424.81 &  387.63 &   42 \\
    SCIPsbt    &  178 &   85.15 &   59.78 &  178 &   42 &  169.78 &  161.44 &   42 \\
    SCIPabt    &  178 &   83.79 &   59.25 &  178 &   42 &  163.36 &  154.26 &   42 \\
    GRBbasic   &  178 &    6.28 &    5.43 &  178 &   42 &   13.03 &   12.04 &   42 \\
    GRBsbt     &  178 &   11.40 &    9.32 &  178 &   42 &   24.88 &   23.29 &   42 \\
    \midrule
    \multicolumn{9}{l}{global budget: $\lceil 0.02 \cdot |E| \rceil$}\\
    SCIPbasic  &  182 &  861.95 &  178.32 &  171 &   20 &  338.87 &  304.96 &   20 \\
    SCIPsbt    &  182 &  340.02 &   78.95 &  178 &   20 &  136.87 &  129.74 &   20 \\
    SCIPabt    &  182 &  373.47 &   79.59 &  178 &   20 &  125.34 &  119.02 &   20 \\
    GRBbasic   &  182 &  692.50 &   25.82 &  175 &   20 &   10.48 &    9.91 &   20 \\
    GRBsbt     &  182 &  292.40 &   22.45 &  179 &   20 &   21.10 &   19.82 &   20 \\
    \midrule
    \multicolumn{9}{l}{global budget: $\lceil 0.03 \cdot |E| \rceil$}\\
    SCIPbasic  &  186 & 1664.22 &  438.26 &  176 &    1 &   67.18 &   67.18 &    1 \\
    SCIPsbt    &  186 &  445.62 &  109.95 &  182 &    1 &   35.41 &   35.41 &    1 \\
    SCIPabt    &  186 &  469.51 &  109.72 &  182 &    1 &   35.42 &   35.42 &    1 \\
    GRBbasic   &  186 & 1246.43 &  111.16 &  179 &    1 &    2.32 &    2.32 &    1 \\
    GRBsbt     &  186 &  365.84 &   57.23 &  183 &    1 &    3.77 &    3.77 &    1 \\
    \midrule
    \multicolumn{9}{l}{global budget: $\lceil 0.04 \cdot |E| \rceil$}\\
    SCIPbasic  &  187 & 1910.56 &  613.64 &  172 &    1 &   63.21 &   63.21 &    1 \\
    SCIPsbt    &  187 &  766.90 &  144.99 &  178 &    1 &   33.13 &   33.13 &    1 \\
    SCIPabt    &  187 &  783.38 &  144.70 &  178 &    1 &   34.23 &   34.23 &    1 \\
    GRBbasic   &  187 & 1457.39 &  114.36 &  173 &    1 &    2.31 &    2.31 &    1 \\
    GRBsbt     &  187 &  327.42 &   64.44 &  187 &    1 &    3.74 &    3.74 &    1 \\
    \midrule
    \multicolumn{9}{l}{global budget: $\lceil 0.05 \cdot |E| \rceil$}\\
    SCIPbasic  &  185 & 1889.19 &  632.42 &  166 &    0 & ---     & ---     & ---  \\
    SCIPsbt    &  185 & 1167.71 &  244.49 &  174 &    0 & ---     & ---     & ---  \\
    SCIPabt    &  185 & 1156.75 &  242.43 &  174 &    0 & ---     & ---     & ---  \\
    GRBbasic   &  185 & 1553.69 &  113.65 &  171 &    0 & ---     & ---     & ---  \\
    GRBsbt     &  185 &  493.57 &   97.94 &  185 &    0 & ---     & ---     & ---  \\
    \midrule
    \multicolumn{9}{l}{global budget: $\lceil 0.06 \cdot |E| \rceil$}\\
    SCIPbasic  &  186 & 1990.98 &  693.58 &  163 &    0 & ---     & ---     & ---  \\
    SCIPsbt    &  186 & 1440.55 &  279.50 &  172 &    0 & ---     & ---     & ---  \\
    SCIPabt    &  186 & 1461.07 &  286.93 &  171 &    0 & ---     & ---     & ---  \\
    GRBbasic   &  186 & 1410.53 &   95.52 &  176 &    0 & ---     & ---     & ---  \\
    GRBsbt     &  186 &  524.40 &  116.61 &  186 &    0 & ---     & ---     & ---  \\
    \midrule
    \multicolumn{9}{l}{global budget: $\lceil 0.07 \cdot |E| \rceil$}\\
    SCIPbasic  &  183 & 2151.99 &  693.12 &  156 &    0 & ---     & ---     & ---  \\
    SCIPsbt    &  183 & 1403.40 &  246.99 &  171 &    0 & ---     & ---     & ---  \\
    SCIPabt    &  183 & 1452.16 &  254.92 &  169 &    0 & ---     & ---     & ---  \\
    GRBbasic   &  183 & 1303.05 &   82.57 &  176 &    0 & ---     & ---     & ---  \\
    GRBsbt     &  183 &  466.03 &  124.70 &  183 &    0 & ---     & ---     & ---  \\
    \midrule
    \multicolumn{9}{l}{global budget: $\lceil 0.08 \cdot |E| \rceil$}\\
    SCIPbasic  &  186 & 2252.93 &  732.66 &  153 &    0 & ---     & ---     & ---  \\
    SCIPsbt    &  186 & 1525.19 &  314.64 &  170 &    0 & ---     & ---     & ---  \\
    SCIPabt    &  186 & 1513.03 &  309.31 &  170 &    0 & ---     & ---     & ---  \\
    GRBbasic   &  186 & 1218.15 &   65.01 &  177 &    0 & ---     & ---     & ---  \\
    GRBsbt     &  186 &  535.35 &  181.75 &  185 &    0 & ---     & ---     & ---  \\
    \midrule
    \multicolumn{9}{l}{global budget: $\lceil 0.09 \cdot |E| \rceil$}\\
    SCIPbasic  &  186 & 2214.69 &  692.20 &  155 &    0 & ---     & ---     & ---  \\
    SCIPsbt    &  186 & 1599.46 &  285.69 &  163 &    0 & ---     & ---     & ---  \\
    SCIPabt    &  186 & 1333.37 &  260.48 &  175 &    0 & ---     & ---     & ---  \\
    GRBbasic   &  186 & 1037.36 &   41.38 &  178 &    0 & ---     & ---     & ---  \\
    GRBsbt     &  186 &  595.20 &  205.80 &  186 &    0 & ---     & ---     & ---  \\
    \midrule
    \multicolumn{9}{l}{global budget: $\lceil 0.10 \cdot |E| \rceil$}\\
    SCIPbasic  &  185 & 2251.07 &  734.90 &  154 &    0 & ---     & ---     & ---  \\
    SCIPsbt    &  185 & 1579.84 &  272.52 &  168 &    0 & ---     & ---     & ---  \\
    SCIPabt    &  185 & 1380.72 &  254.21 &  173 &    0 & ---     & ---     & ---  \\
    GRBbasic   &  185 & 1153.80 &   41.88 &  175 &    0 & ---     & ---     & ---  \\
    GRBsbt     &  185 &  578.28 &  217.22 &  185 &    0 & ---     & ---     & ---  \\
    \bottomrule
  \end{tabular*}
  \end{scriptsize}
  \end{center}
  \vskip -0.1in
\end{table}


\end{document}